\documentclass[a4paper]{article}
\usepackage[utf8]{inputenc}
\usepackage[intlimits]{amsmath}
\usepackage{amssymb, units, amsthm, bm, comment, esint, xcolor}
\usepackage[T1]{fontenc}
\usepackage{wrapfig}
\usepackage{hyperref}
\usepackage{srcltx,
}
\usepackage{graphicx}

\newtheorem{thm}{Theorem}[section]

 \newtheorem{prop}[thm]{Proposition}
 \theoremstyle{definition}
 \newtheorem{defn}[thm]{Definition}
 \theoremstyle{remark}
 \newtheorem{rem}[thm]{Remark}
 
 \numberwithin{equation}{section}
\def\Pr{\noindent {\bf Proof. }}
\def\diver{\mathop{\mathrm{div}}\nolimits}
\def\tens#1{\mathbb{#1}}
\def\vec#1{\boldsymbol{#1}}
\def \R{\mathbb{R}}
\def \N{\mathbb{N}}
\def \O{\Omega}

\def\chiu{\hfill$\displaystyle\vspace{4pt}
\underset\Box\null$\par}

\def\bT{\tens{T}}
\def\bD{\tens{D}}
\def\bG{\tens{G}}
\def\bO{\tens{O}}
\def\bS{\tens{S}}
\def\bI{\tens{I}}
\def\bZ{\tens{Z}}

\def\bW{W}
\def\bB{\tens{B}}
\def\bA{\tens{A}}

\def\bx{\vec{x}}
\def\bu{\vec{u}}
\def\bv{\vec{v}}
\def\bw{\vec{w}}
\def\bs{\vec{s}}

\def\bn{\vec{n}}
\def\bb{\vec{b}}
\def\bz{\vec{z}}
\def\bnul{\vec{0}}

\def\bh{\vec{h}}

\def\pa{\partial}
\def\na{\nabla}
\def\be{\begin{equation}}
\def\ba{\begin{array}}
\def\ea{\end{array}}
\def\ee{\end{equation}}
\def\displ{\displaystyle\vspace{6pt}}
\def\btau{\vec{\tau}}

\def\Lnd#1{L^{#1}_{\bn, \diver}}

\def\Wn#1{W^{1,#1}_{\bn}}
\def\Wnd#1{W^{1,#1}_{\bn, \diver}}
\def\Wndm#1{W^{-1,#1}_{\bn, \diver}}
\def\Wnm#1{W^{-1,#1}_{\bn}}

\def\ss{_\mathrm{s}}
\def\ff{_\mathrm{f}}

\def\rhotf{\rho^\mathrm{m}\ff}
\def\rhots{\rho^\mathrm{m}\ss}

\def\vf{\bv\ff}
\def\vs{\bv\ss}
\def\pa{\partial}
\def\Tf{\mathbb{T}\ff}
\def\Ts{\mathbb{T}\ss}
\def\Ttot{\mathbb{T}}
\def\Stot{\mathbb{S}}
\def\pf{p_{\rm f}}
\def\II{\mathrm{\bf I}}
\def\pft{p\ff^\mathrm{t}}

\def\pef{p^\mathrm{eff}}
\def\muf{\mu\ff}

\begin{document}
\title{On three-dimensional flows of pore pressure activated Bingham fluids}
\author{ A.~Abbatiello\thanks{
Dipartimento di Matematica e Fisica, Universit\`{a} degli Studi della Campania ``Luigi Vanvitelli'', viale Lincoln n.5, 81100 Caserta,
 Italy. {\it E-mail address}: anna.abbatiello@unicampania.it,}, T.~Los\thanks{Charles University, Faculty of Mathematics and Physics, Mathematical Institute, 186
75 Prague 8, Czech Republic. {\it E-mail addresses}:  los@karlin.mff.cuni.cz, malek@karlin.mff.cuni.cz, soucek@karel.troja.mff.cuni.cz.}, J.~M\'{a}lek\footnotemark[2], O.~Sou\v{c}ek\footnotemark[2],}
\date{}
\maketitle
\textit{To Professor Vsevolod Alekseevich Solonnikov on the occasion of his 85th birthday}
\begin{abstract}
We are concerned with a system of partial differential equations describing internal flows of homogeneous incompressible fluids of Bingham type in which the value of activation (the so-called yield) stress depends on the internal pore pressure governed by an advection-diffusion equation. 
After providing the physical background of the considered model, paying attention to the assumptions involved in its derivation, we focus on the PDE analysis of the initial and boundary value problems. 
We give several equivalent descriptions for the considered class of fluids of Bingham type. In particular, we exploit the possibility to write such a response   as an implicit tensorial constitutive equation, involving the  pore pressure, the deviatoric part of the Cauchy stress and the velocity gradient. Interestingly, this tensorial response can be characterized by two scalar constraints. We employ a similar approach to treat stick-slip boundary conditions. 
Within such a setting we prove long time and large data existence of weak solutions to the evolutionary problem in three dimensions.
\end{abstract}

\section{Introduction}

The mechanical behavior of water saturated geological materials such as soils or sands is known to involve the notion of the so-called effective stress (or effective pressure), introduced in 1920 by Terzaghi \cite{terzaghi-1925}. The effective pressure is defined as the difference between the mean normal stress in the medium and the pressure of the interstitial fluid - {\it  pore pressure}. Water saturated geological materials are mixtures composed of an unconsolidated granular solid material and an interstitial pore space occupied by a fluid, see  Fig.~\ref{fig:fig1}. With this picture in mind the total stress exerted on any control surface in such a medium comprises two contributions, namely the stress transmitted by the fluid and the stress transmitted { by the granular solid}. Mechanical loading or unloading of such a saturated material due to external forcing leads to redistribution of the stresses between the two { constituents}, which in general can be a rather complex process. Despite its complexity, several general observations can be made. First, if during the process the stress in the granular material increases, for example by reducing the pore pressure while keeping the total loading constant, the granular structure compactifies and becomes more rigid. A textbook example of this process is the beach sandcastle stabilization, when the fluid flowing out of the wet sand stabilizes the sand by ``sticking'' the sand grains closer together (here also capillary phenomena play a significant role).  Second, as an opposite extreme, it may happen that during some processes the pressurized interstitial fluid bears almost the whole mechanical load exerted on the system, which leads to effective mechanical decoupling of the solid grains and the so called ``liquefaction'' can occur.

{ In this paper we develop a mathematical theory for a model that can be viewed as a simple toy model for the process of pore-pressure activated flows of saturated granular materials  decribed above. We give up the ambition to model the actual process of liquefaction of real-world geological materials such as soils, since compared to what is presented here, this would require much more involved modelling of the activation yield criteria for such materials and of their rheological properties after the activation. However, we believe that even the strongly simplified setting presented here provides certain qualitative insight into the physics of pore-pressure activated flows and may even have some relevance to the problems of static liquefaction (see \cite{Lade}) or enhanced oil recovery.\footnote{{ In enhanced oil recovery steam or carbon dioxide is injected to reclaim oil that remains after initial extraction (see \cite{kirken} for a discussion of enhanced oil recovery and \cite{NakRaj}, \cite{sribonraj} for modeling and numerical studies). The recovery takes place after a pressure builds-up in the porous substrate containing the remnant oil and the oil starts to flow. Before the flow takes place we do have the steam/carbon dioxide being pumped into the porous rock and this flow is governed by some Darcy-like equation. The pressures involved are quite high and the material properties like the viscosity of the fluid would be pressure dependent, and at such high pressures the porous rock would undergo some deformation, these two effects are being ignored. We are also not modeling the porous rock as an individual constituents, as the considered mixture is constituted by steam and oil.}}
All this in our view justifies to study the associated initial and boundary value problems in terms of mathematical well-posedness, which represents the main objective of the manuscript. 

The model developed here is obtained within the context of the theory of interacting continua initiated by Truesdell \cite{Truesdellmix}, \cite{Truesdell} 
(see also the review articles by Bowen \cite{Bowen}, Atkin and Craine \cite{atkincraine}, and the numerous appendices in the book on rational thermodynamics by Truesdell \cite{Truesdell_rm}, and the books by Samoh\'{y}l \cite{Sam}, Rajagopal and Tao \cite{RajTao}). Within this framework, we are concerned with the flow of a mixture composed of two fluid components, one representing the unconsolidated granular material flowing once a certain activation critierion is met, and the second fluid being Newtonian, representing the pore-space fluid. The Darcy-type flow of the pore fluid relative to the second fluid is considered, driven by the pore pressure gradient and gravity. This flow accomodates the pore presssure. Once the pore pressure reaches a certain threshold, the ``granular'' fluid starts to flow. 



\par

The organization of the paper is as follows.} In Section 2, we develop the model from the  principles within the framework of mixture theory (see \cite{Truesdell, Bowen, RajTao, Sam, Drew-Passman-1998, MR}) through a number of physically reasonable approximations. In particular, we compare our final system with the one studied recently by Chupin and Math\'e \cite{CM}; they differ by the structure of the right-hand side in the equation for the (fluid/effective) pressure. In Section 3, we reformulate the response as an implicit equation involving the constitutively determined part of the stress, the symmetric part of the velocity gradient and the pore pressure. We also give an alternative characterization of such implicit constitutive equation in terms of two scalar constraints (extending here an interesting observation from \cite{CM}). 
Next, we focus on the mathematical analysis of the initial and boundary value problem in three-dimensional domains. Note that Chupin and Math\'e analyzed only two-dimensional flows in \cite{CM}, which is easier as the energy equality holds for a weak solution.
We consider internal flows when the whole boundary is impermeable and we study the problem with stick-slip boundary conditions (that can be also equivalently written as an implicit constitutive equation on the boundary and characterized by two inequalities). Stick-slip (or threshold slip) states that the velocity does not slip until the amplitude of the tangent part of the normal traction on the boundary  exceeds a certain critical value. This boundary condition, which is physically relevant to the pore pressure activated fluids considered in the bulk, includes Navier's slip and (perfect) slip boundary conditions as special cases. We establish the long-time and large-data existence of the corresponding weak solutions; see Section 4 for the formulation of the main result and Section 6 for its proof. We exploit the characterization of the implicit constitutive equations by two scalar constraints,  both in the bulk and on the boundary, as a  tool to show that the limit object of suitable approximative sequences fulfils these constitutive equations as well, see Proposition \ref{convergence-lemma} proved in Section 5, where we also introduce the approximations and study their properties. Finally, we comment on possible results for no-slip boundary conditions and further extensions in the concluding section. 

\section{Derivation of the model in the framework of multi-constituent theory}

\begin{wrapfigure}{h!}{0.435\textwidth}
\includegraphics[scale=0.065]{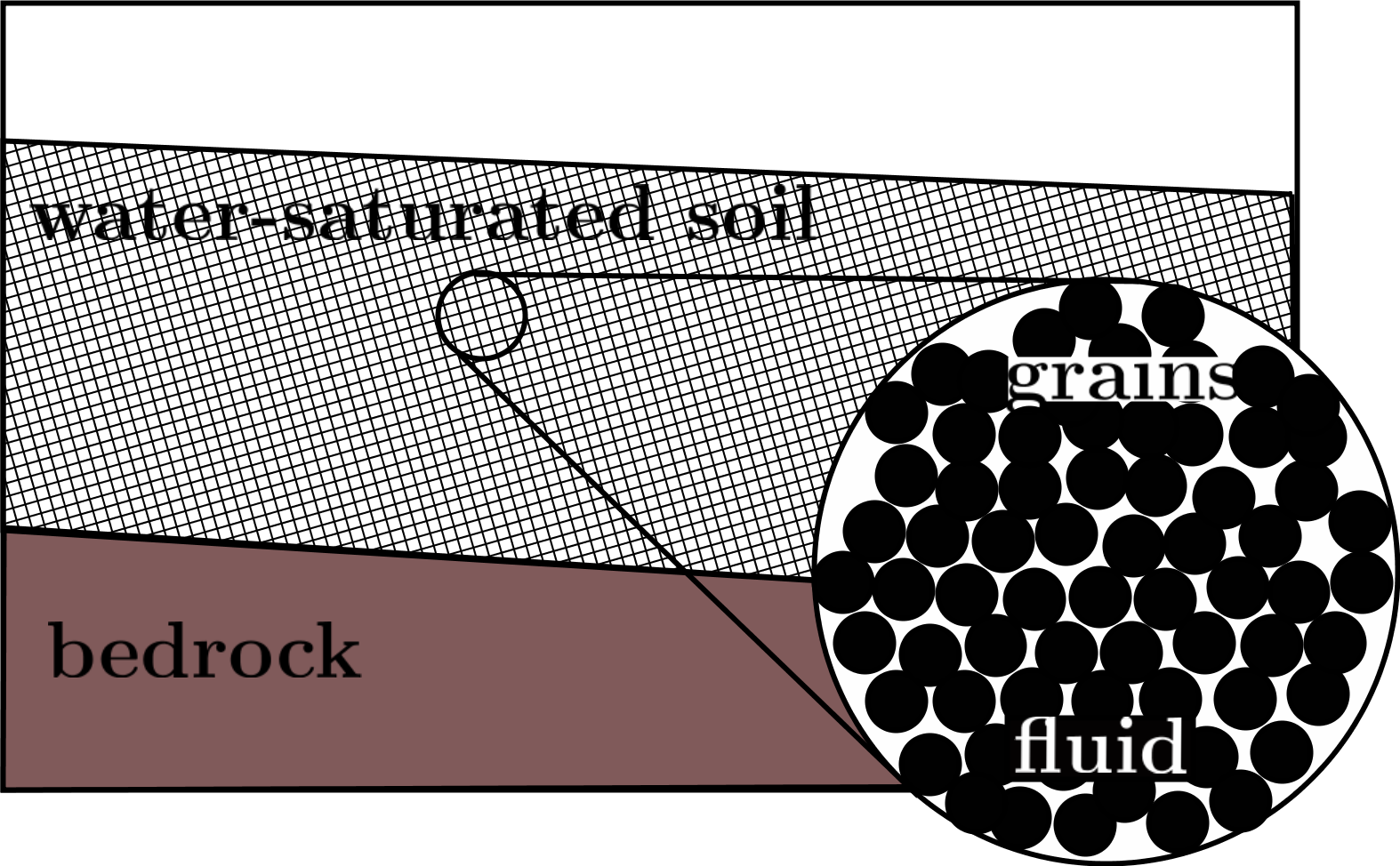}
\caption{\footnotesize{Sketch of the typical problem geometry and zoom into the structure of the material composed of { a graular unconsolidated solid} filled with an interstitial fluid.}}
\label{fig:fig1}
\end{wrapfigure}
We consider a { continuum} that is  composed, 
in a homogenized sense of two types of materials: a { granular unconsolidated solid (fllowing once the stress exceeds the value of activation stress)} and the interstitial pore space occupied  by a fluid, see in Fig.~\ref{fig:fig1}. These two materials are called { constituents}.  { We start with general description of such a mixture and then specify a number of simplifying physical assumptions that will result to a model we analyze in the remaining part of the paper. { We refer to quantities related to the
granular material (flowing after an activation criterion takes place) as {\it solid} (denoted by subscript "$\mathrm{s}$") and to the interstitial fluid simply as {\it fluid} (denoted by subscript "$\mathrm{f}$").}}

Based on the theory of multi-component materials (see e.g. \cite{RajTao} or \cite{MR}), we first formulate the individual mass and momentum balances for both components. Restraining ourselves to a purely mechanical setting, for simplicity, we do not need to formulate the balance equations for energy and entropy. 
\par
The  balance equations for mass read as follows:
\begin{subequations}
\label{mass-balances}
\begin{align}
\label{mass-balance-f}
	\frac{\pa(\phi\rhotf)}{\pa t} + \diver({\phi\rhotf\vf}) &= 0\ ,\\
\label{mass-balance-s}
	\frac{\pa((1{-}\phi)\rhots)}{\pa t} + \diver((1{-}\phi){\rhots\vs}) &= 0\ .
\end{align}
\end{subequations}
Here $\rhotf$ and $\rhots$ denote the material (true) densities of the fluid and the solid, $\phi$ denotes the volume fraction of the {fluid} (equal to { the porosity of the granular solid} in the saturated case considered here) and $\vf$ and $\vs$ denote the velocities of the { constituents}, respectively. The zero on the right-hand side of equations \eqref{mass-balances} expresses the fact that we do not consider any mass transfer between the { constituents.} 
\par The balance equations for  linear momentum for these two { constituents} take the form
\begin{subequations}
\label{momentum-balances}
\begin{align}
	\frac{\pa(\phi\rhotf\vf)}{\pa t} + \diver({\phi\rhotf\vf{\otimes}\vf}) &= \diver\Tf + \phi\rhotf\bb + \II\ ,\\
	\frac{\pa((1{-}\phi)\rhots\vs)}{\pa t} + \diver({(1{-}\phi)\rhots\vs{\otimes}\vs}) &= \diver\Ts + (1{-}\phi)\rhots\bb - \II\ ,
\end{align}
\end{subequations}
where, for $\bu,\bw\in \mathbb{R}^3$ the symbol $\bu\otimes\bw$  denotes the tensor of components $(\bu\otimes\bw)_{ij}:=u_iw_j$ with $i,j=1, 2, 3$, while $\Tf$ and $\Ts$ stand for the fluid and solid Cauchy stresses, respectively, both of which are assumed to be symmetric (i.e. $\Tf{=}\Tf^\mathrm{T}$, $\Ts{=}\Ts^\mathrm{T}$). The quantity $\mathrm{\II}$ represents the interaction force between the { constituents}. The interaction nature of this force is reflected by the fact that it appears with plus sign in one equation and with minus in the other. Finally $\bb$ is the body force (same for both { constituents}, typically this is the gravity acceleration vector).
In reality, both { constituents} are compressible, i.e. both material densities $\rhots$ and $\rhotf$ must be specified by a corresponding state equation. In the isothermal setting considered here, such relation would take the form of dependence on the material stress state of the particular { constituent}. Since the dominant compressibility effect in the context of real-world  geological materials is not related to the changes of material densities, but rather to the changes in porosity in reaction to the applied loading (see \cite{Domenico-Schwartz-1998}, chapter 4), we neglect the former effect by setting 
\begin{align}
	\label{constant-rhot}
	\rhotf = \mathrm{const}\ff \ \ \  \mbox{and} \quad \rhots=\mathrm{const}\ss.
\end{align}
Dividing now \eqref{mass-balance-f} by $\rhotf$ and \eqref{mass-balance-s} by $\rhots$ and summing the resulting equations, we obtain
\begin{subequations}
\label{mass-bal-2}
\begin{align}
	\label{mb-diver}
	\diver\vs = -\diver(\phi(\vf{-}\vs))\ .
\end{align}
Inserting this relation into \eqref{mass-balance-s} (divided by $\rhots$) yields the evolution equation for porosity
\begin{align}
	\label{mb-phi}
	\frac{\pa\phi}{\pa t} +  \vs\cdot\nabla\phi = -(1{-}\phi)\diver\left({\phi(\vf{-}\vs)}\right)\ .
\end{align}
\end{subequations}
Under the assumptions \eqref{constant-rhot}, the system of equations \eqref{mass-bal-2} is equivalent to the system \eqref{mass-balances}.
\par
The balance equations for linear momentum \eqref{momentum-balances} are reformulated in terms of an equivalent system, where the balance equation for linear momentum of the solid  is replaced by the balance equation for linear momentum of the mixture as a whole. Thus, using \eqref{constant-rhot}, we get
\begin{subequations}
\label{momentum-balances-22}
\begin{align}
\label{mom-bal-fluid}
	\rhotf\left(\frac{\pa(\phi\vf)}{\pa t} + \diver({\phi\vf{\otimes}\vf})\right) &= \diver\Tf + \phi\rhotf\bb + \II\ ,\\\nonumber
	\rhots\left(\frac{\pa\vs}{\pa t} + \diver{(\vs{\otimes}\vs)}\right) &= \diver\Ttot + \rho\bb  - \frac{\pa}{\pa t}\left(\phi(\rhotf\vf{-}\rhots\vs)\right)\label{mom-bal-tot1}\\ & - \diver(\phi(\rhotf\vf{\otimes}\vf -\rhots\vs{\otimes}\vs))\ ,
\end{align}
\end{subequations}
where, in the second equation, we introduced the total Cauchy stress $\Ttot$ and the total density $\rho$ by
\begin{align}\label{density-ms}
	\Ttot := \Ts + \Tf, \hspace{1cm} \rho:= \phi\rhotf + (1{-}\phi)\rhots\ =\rhots+\phi(\rhotf-\rhots).
\end{align} 

Next we introduce several simplifications.
\begin{itemize}
\item 
	In the balance of linear momentum for the fluid \eqref{mom-bal-fluid}, we ignore the inertial forces, i.e. the whole left-hand side of \eqref{mom-bal-fluid} is set to zero. We further consider 
	$\Tf$ of the form
	\begin{align}
		\label{fluid-Cauchy-simplification}
		\Tf =  -\pft\phi\,\mathbb{I}\ ,
	\end{align}
	where $\pft$ is the true pressure in the interstitial fluid (pore pressure).
	Finally, the interaction force $\II$ takes a simple form corresponding to the linear drag
	\begin{align}
		\label{interaction-force}
		 \II = -\alpha(\vf{-}\vs) + \pft\nabla\phi\ ,
	\end{align} 
	where $\alpha$ is the drag coefficient of the form
	\begin{align}
\label{alpha-darcy}
	\alpha:=\frac{\phi^2\muf}{k(\phi)}\ ,
\end{align}
 $\muf$ being the dynamic viscosity of the fluid { (assumed to be constant for simplicity)} and $k(\phi)$  the permeability of the granular material. The presence of the second term on the right-hand side of \eqref{interaction-force}  is known from multiphase continuum theory as an artefact of the volume averaging technique \cite{Drew-Passman-1998}, which must be present to cancel out in the fluid momentum balance with a corresponding term coming from the divergence of eq. \eqref{fluid-Cauchy-simplification}. See also \cite{MR} where such terms occur from the derivation directly.
	
\item  In the balance equation \eqref{mom-bal-tot1}, we keep the inertial term only on the left-hand side and  neglect the last two terms on the right hand side using a rough scaling argument stating that the scale of these terms is at most the scale of the left hand side, multiplied by the scale of porosity, which, in the considered applications, typically does not exceed a few percent. Furthermore, since $\phi$ is typically below $0.1$, we conclude that 
 $\rho$ { introduced in \eqref{density-ms}} is approximately equal to $\rhots$, so { we replace $\rho\bb$ by $ \rhots\bb$ in \eqref{mom-bal-tot1}}. 
\end{itemize}
With these  simplifications, the balance equations \eqref{mass-bal-2} and  \eqref{momentum-balances-22}   take the form
\begin{subequations}\label{210}
\begin{align}
\frac{\partial \phi}{\partial t}+\vs\cdot\nabla\phi&=-(1-\phi)\diver\left(\phi(\vf-\vs)\right),\label{a}\\
\diver\vs&=-\diver(\phi(\vf-\vs)),\label{b}\\
\phi\nabla\pft&=\phi\rhotf\bb-\alpha (\vf-\vs),\label{c}\\
\rhots\left(\frac{\pa\vs}{\pa t} +\diver(\vs\otimes\vs)\right)&=\diver\Ttot+{\rhots}\bb.
\end{align}
\end{subequations}
 Next, we  multiply \eqref{c} by $\frac{\phi}{\alpha}$, use \eqref{alpha-darcy} and apply divergence to the result. After inserting the outcome of these computations in  \eqref{b} we obtain
\begin{equation}\label{newb}
\diver\vs = \diver\left(\frac{k(\phi)}{\muf}(\nabla p_f^t - \rhotf\bb)\right).
\end{equation}
As a consequence,  \eqref{newb} replaces \eqref{b}. 
\par
The next assumption 
states that the porosity of the granular solid can be described by a constitutive relation of the form
\begin{align}
\label{porosity-effective-stress}
\phi=\widehat{\phi}(\pef)\, \hspace{1cm}\mathrm{where}\ \hspace{1cm}\pef := p-\pft.
\end{align}
The quantity $\pef$\!, called  {\it effective pressure} as introduced by Terzaghi \cite{terzaghi-1925},  is defined as the difference between the total mixture pressure and the fluid (pore) pressure. The quantity $\pef$ is assumed to reflect the part of the loading bore by the granular solid. Inserting  the constitutive assumption \eqref{porosity-effective-stress} in \eqref{a} and using { \eqref{b} and} \eqref{newb}  we obtain the following evolution equation for the effective pressure
\begin{align}
\label{eq-p-eff-evol}
\frac{d\widehat{\phi}}{d\pef}\left(\frac{\pa\pef}{\pa t} + \vs\cdot\nabla\pef\right) = (1{-}\phi)\diver\left(\frac{k(\phi)}{\muf}(\nabla\pft - \rhotf\bb)\right)\!\!.
\end{align}
Setting 
\begin{align}
-\frac{1}{\beta}:=\frac{d\widehat{\phi}}{d\pef} \ \ \mbox{with}  \ \beta>0,
\end{align}
replacing $\pft$ in \eqref{newb} and \eqref{eq-p-eff-evol} by $p-\pef$ and 
splitting the total Cauchy stress $\mathbb{T}$ as
\begin{align}
\Ttot = -p\mathbb{I} + \Stot,
\end{align}
we arrive at the following set of governing equations
\begin{subequations}
\label{system-I}
\begin{align}
\label{sysI-divers}
\diver\vs &= \diver\left(\frac{k(\phi)}{\muf}(\nabla p-\nabla \pef - \rhotf\bb)\right)\ ,\\
\label{sysI-momeq}
\rhots\left(\frac{\pa\vs}{\pa t} + \diver{(\vs{\otimes}\vs)}\right) &= -\nabla p +\diver\Stot + \rhots\bb\ ,\\\label{sysI-phieq}
\frac{\pa\pef}{\pa t} + \vs\cdot\nabla\pef &=- (1{-}\phi)\beta\diver\left(\frac{k(\phi)}{\muf}(\nabla p-\nabla\pef - \rhotf\bb)\right)\!, \\
\label{sysI-vf} { \vf} &{= \vs - \frac{1}{\alpha}\widehat{\phi}(\pef){\left(\nabla\pft - \rhotf\bb \right)\,.}}
\end{align}
\end{subequations}
{ Since $\pft = p- \pef$, we can view \eqref{system-I} as the system of partial differenctial equations describing the evolution of $p$, $\vs$, $\pef$ and $\vf$, where t}he rheology of the material needs to be specified by providing a constitutive relation for the stress $\bS$. We shall assume that the material behaves as very stiff until the threshold is reached at which moment the { material} starts to flow as a liquid. A simplest constitutive relation for this type of response is characterized by the so-called Bingham fluid \cite{Bingham:1916ki}, where the solid part responses as a perfectly rigid body until the magnitude of the stress exceeds the threshold when the solid flows as a Newtonian fluid. This type of response is usually written in the following way: 
\be \label{def-sigma}
\begin{cases}\vspace{6pt}\displaystyle
|\bS|\leq \tau(\pef) \ &\mbox{if and only if} \  \bD=\bO ,\\ \vspace{6pt}\displaystyle
|\bS|> \tau(\pef)   &\mbox{if and only if} \  \bS =\tau(\pef)\frac{ \bD}{|\bD|}+2 \nu_* \bD.
\end{cases}\ee
Here $\bD$ is the symmetric part of the velocity gradient 
$$\bD:=\frac{1}{2}(\nabla\vs+(\nabla\vs)^\mathrm{T})\ ,$$
$\nu_*>0$ is the viscosity and $\tau({\pef})$ is the threshold, depending on the effective pressure.
Typically,  
\begin{equation}\label{sysI-tau}
	\tau({\pef}) = q_0(\pef)^+,
\end{equation}
where $q_0$ is a constant and the symbol $( )^+$ denotes the positive part of a quantity, i.e. $(\psi)^+ := \mathrm{max}(\psi, 0)$.

{ The activation criterion \eqref{def-sigma} is too simple to describe the shear instability of real-world granular materials since it does not take into account any concept of internal friction. In reality, it should be replaced by some form of Mohr-Coulomb criteria, see, e.g., \cite{Atkinson}.
 Similarly, also} {the fact that the material, which is a mixture of flowing solid particles and fluid (a slurry), is supposed to response, after being activated,  as a Navier-Stokes (linear viscous) fluid, is a severe limitation of the model considered here. To overcome this defficiency one would need to incorporate more realistic models used for description of flows of granular material (capable of exhibiting normal stress differences, etc.). Such models have been developed in \cite{johnson1}, \cite{johnson2}, see also a review article \cite{hutter}.}


{ Despite these important limitations,} the system \eqref{system-I} -- \eqref{sysI-tau} seems to be a meaningful, physically justified and relatively simple model worth of studying. We however do not further investigate this system here,  as our goal  is to identify the assumptions that lead to the model analyzed in \cite{CM}. Towards this aim, we introduce the following additional assumptions:
\begin{itemize}

\item{\it Pore pressure evolution approximation.}  Using the relation $\pft= p-\pef$ we rewrite \eqref{sysI-phieq} as an evolution equation for $\pft$ 
\begin{align}
\label{eq-p-fluid-evol}
\frac{\pa\pft}{\pa t} + \vs\cdot\nabla\pft = \frac{\pa p}{\pa t} + \vs\cdot\nabla p+(1{-}\phi)\beta\diver\left(\frac{k(\phi)}{\muf}(\nabla\pft - \rhotf\bb)\right).
\end{align}
Next we assume that the dominant contribution to the total pressure $p$ in eq.~\eqref{eq-p-fluid-evol} comes from the hydrostatic part, which may in general depend explicitly on time to include problems with evolving boundary. Consequently, we replace $p({x}, t)$ by $p_s( x, t)$ in eq.~\eqref{eq-p-fluid-evol}, where $p_s$ is a given function. Also, we replace $(1{-}\phi)$ by $1$ on the right-hand side of \eqref{eq-p-fluid-evol} since, as set above, we are interested in situations where $\phi<0.1$ and finally, we assume that both the permeability $k$ and the compressibility parameter $\beta$ are constant. Setting thus 
\begin{align}
	K := \frac{\beta k}{\muf} \geq 0\ ,
\end{align}
the equation  \eqref{eq-p-fluid-evol} simplifies to the form
\begin{align}
\frac{\pa\pft}{\pa t} + \vs\cdot\nabla\pft &= K\Delta\pft - \diver(K\rhotf\bb) + {\frac{\pa p_s}{\pa t}} + \vs\cdot\nabla p_s.
\end{align}
\item {\it Yield criterion approximation.} Also in the yield criterion, we replace the pressure $p$ in the definition of the effective pressure \eqref{porosity-effective-stress} by $p_s$, i.e. instead of  \eqref{sysI-tau}, we have
\begin{align}
	\tau({\pft}) = q_0(p_s - \pft)^+\ .
\end{align}

\item {\it Incompressibility.} We ignore the effect of porosity changes in  \eqref{sysI-divers}  by replacing \eqref{sysI-divers}  with the incompressibility constraint
\begin{align}
\diver\vs = 0\ .
\end{align}
\end{itemize}
With the above set of simplifying assumptions, the final reduced system of governing equations reads as follows
\begin{subequations}
\label{system-II}
\begin{align}
\label{sysII-divers}
\diver\vs &= 0\ ,\\
\label{sysII-momeq}
\rhots\left(\frac{\pa\vs}{\pa t} + \diver{(\vs{\otimes}\vs)}\right) &= -\nabla p +\diver\Stot + \rhots\bb\ ,\\
\label{sysII-phieq}
\frac{\pa\pft}{\pa t} + \vs\cdot\nabla\pft &= K\Delta\pft  - \diver(K\rhotf\bb) + {\frac{\pa p_s}{\pa t}} + \vs\cdot\nabla p_s\ ,
\end{align}
\end{subequations}
where $\bS$ satisfies
\be \label{constitutive-model-II}
\begin{cases}\vspace{6pt}\displaystyle
|\bS|\leq \tau(\pft) \ &\mbox{if and only if} \  \bD=\bO ,\\ \vspace{6pt}\displaystyle
|\bS|> \tau(\pft)   &\mbox{if and only if} \  \bS =\tau(\pft)\frac{ \bD}{|\bD|}+2 \nu_* \bD,
\end{cases} 
\ \ \mbox{with } \ \tau({\pft}) = q_*(p_s- \pft)^+,
\ee
{ and where the velocity $\vf$ is given by} 
\be \label{system-II-+}
{ \vf = \vs - \frac{1}{\alpha}\widehat{\phi}(p-\pft)\left(\nabla\pft - \rhotf\bb \right)\,.}
\ee
{ Since $\vf$ does not enter into \eqref{system-II} and \eqref{constitutive-model-II}, the equation \eqref{system-II-+} describing the evolution of $\vf$ is not considered anymore in what follows (as $\vf$ can be always obtained from equation of Darcy's type \eqref{system-II-+} once $\vs$ and $\pft$ are known/computed from \eqref{system-II} and \eqref{constitutive-model-II}.}
 
\section{(Re)-formulations of the Problem}
Let  $T$ be a positive real number and $\O\subset \R^3$  a bounded domain with the boundary $\pa\O$. We set $Q_T:=(0,T)\times \O$ and $\Sigma_T:= (0,T)\times\pa\O$. The symbol  $\bn :\pa\O\to \mathbb{R}^3$ denotes the outer unit normal vector while  for any vector $\bz$ defined on $\pa\O$ we set $\bz_\tau:=\bz-(\bz\cdot \bn)\bn$ representing the projection of $\bz$ to the tangent plane.
 \par We consider  unsteady flows of a homogeneous incompressible non-Newtonian fluid of a Bingham type with a variable threshold,  described in the previous section, see \eqref{system-II}-\eqref{constitutive-model-II}. In what follows, we slightly change the notation and write $\bv$ instead of $\vs$, $\varrho_*$ instead of $\rhots$ and $\pf$ instead of $\pft$. We also set $g:= {\frac{\pa p_s}{\pa t}} - \diver(K\rhotf\bb)$.
\par
Following a recent observation in \cite{BGMS} (see also \cite{RS}, \cite{BM}, \cite{BM17}) the rheological behaviour \eqref{constitutive-model-II}  can be equivalently written as
\be\label{model1}
2 \nu_*\bD =  \frac{\left(|\bS| - \tau(\pf)\right)^{+}}{|\bS|} \bS \ \ \mbox{where}\ \ \ \tau({\pf}) = q_*(p_s- \pf)^+.
\ee
 
We are thus interested in solving the following problem. For given $\varrho_*, \nu_*,q_*\in (0, \infty)$, $\bb:Q_T \rightarrow \R^3$, $g:Q_T\to \R$, $p_s: Q_T\to\R$, we look for $\bv:Q_T  \rightarrow \R^3$, $p, \pf :Q_T  \rightarrow \R$ and $\bS:Q_T  \rightarrow \R^{3\times 3}_{\textrm{sym}}$ satisfying 
\be\label{system}
\ba{l}\displ
\diver\bv = 0 \textrm{ in } Q_T, \\ \displ
\varrho_*\left(\pa_t\bv + \diver(\bv\otimes \bv)\right) - \diver \bS  +\nabla p = \varrho_* \bb \textrm{ in } Q_T, \\\displ
\bG(\bS, \bD, \pf) = \bO \textrm{ in } Q_T, \\\displ
\pa_t \pf+ \bv \cdot \na \pf- K \Delta \pf=g+\bv\cdot\nabla p_s \textrm{ in } Q_T, 
\ea
\ee
where $\bG$ is a continuous function defined on $\R_{\rm sym}^{3\times3}\times\R_{\rm sym}^{3\times3}\times\R$ through 
\begin{equation}
\bG(\bS, \bD, \pf) =  \frac{\left(|\bS| - q_*(p_s- \pf)^+\right)^{+}}{|\bS|} \bS-2 \nu_*\bD.
\end{equation}
In addition, the unknown functions $(\bv, p, \pf, \bS)$ are required, for given $\sigma_*, \gamma_*\in [0, \infty)$ and $\bv_0:\Omega \to\R^3$ (such that $\diver\bv_0=0$ in $\Omega$) and $p_0:\Omega\to\R$, to fulfil the following initial and boundary conditions: 

\begin{alignat}{1}
&\bv(0, \cdot)=\bv_0 \ \mbox{and} \  \pf(0,\cdot)={p}_0 \ \mbox{in} \ \O,\\ 
\label{boundary}
&\bv\cdot\bn=\bnul \ \ \mbox{and} \ \ \na \pf \cdot \bn=\bnul   \  \ \mbox{on} \ \ \Sigma_T,\\
\label{g2}
&\gamma_*\bv_\tau=\frac{(|\bs|-s_*)^+}{|\bs|}\bs \ \ \mbox{on} \ \Sigma_T,
\end{alignat}
where $s_8\in (0, \infty)$.
The conditions in \eqref{boundary} state that the boundary is impermeable, while \eqref{g2} characterizes the result of the interaction of the fluid and the boundary along the boundary. Here 
$$ \bs:=-(\bT\bn)_\tau=-(\bS\bn)_\tau.$$
Note that  \eqref{g2} usually written as
\begin{equation}\label{delta1}
\begin{cases}
|\bs|\leq s_* \ &\mbox{if and only if} \ \bv_\tau=\bnul,\\
|\bs| { >} s_* \ &\mbox{if and only if} \ \bs=s_*\frac{\bv_\tau}{|\bv_\tau|}+\gamma_*\bv_\tau,
\end{cases}
\end{equation}
 describes the stick-slip (or threshold slip) and includes, as special cases,  Navier's slip condition by taking $s_*=0$ and $\gamma_*>0$, and perfect slip condition if in addition $\gamma_*=0$. Note that the no-slip condition is obtained by letting either $s_*\to +\infty$ or $\gamma_*\to +\infty$.
\par
One of the motivations for this work is a recent paper by Chupin and Math\'e \cite{CM} where the authors characterize the tensorial response \eqref{model1} through two scalar constraints:
\begin{equation}\label{model2-2}
|\bZ|\leq \tau(\pf) \ \ \mbox{and} \  \ \bZ :\bD\geq \tau(\pf) |\bD| \ \ \mbox{where} \ \  \bZ:= \bS-2 \nu_* \bD \ \mbox{and} \ \tau(\pf)=q_*(p_s- \pf)^+\!.\end{equation}

In fact, Chupin and Math\'e \cite{CM} considered the second constraint with the equality sign in their existence result concerning planar flows, but then they incorrectly argue when performing the limit in the constitutive equation (see Step 2 (a) in \cite{CM}). This difficulty can be overcome easily if the inequality is used here instead of the equality, as shown in the proof of Theorem \ref{main-thm} in Section \ref{sectionProof} below.

\par Before we prove that \eqref{model1} and \eqref{model2-2} are equivalent, we provide analogously a condition that characterizes \eqref{g2}. It takes the form 
\begin{equation}\label{model2-bc}
|\bz|\leq s_* \ \mbox{and} \ \  \bz\cdot \bv_\tau\geq s_*|\bv_\tau| \ \ \mbox{where} \ \bz:=\bs-\gamma_*\bv_\tau.
\end{equation}
Next, we prove the following statement. 
\begin{prop}\label{equivalent-formulation}The following equivalences hold:
\begin{itemize}
\item[(a)] \eqref{constitutive-model-II} $\iff$ \eqref{model1} $\iff$ \eqref{model2-2};
\item[(b)]  \eqref{delta1}  $\iff$ \eqref{g2} $\iff$ \eqref{model2-bc}.
\end{itemize}
\end{prop}
\Pr The equivalence \eqref{constitutive-model-II} $\iff$ \eqref{model1} is simple. We prove that \eqref{model1} is equivalent to \eqref{model2-2}. Let us first assume that $(\bS, \bD, \pf)$ fulfil \eqref{model1}. If $\bD=\bO$ then $|\bS|\leq \tau(p_f)$ and $\bZ=\bS$ and \eqref{model2-2} holds. If  $\bD\neq \bO$, then $|\bS|>\tau(p_f)$, and the formula \eqref{model1} implies
\begin{equation}\label{ZS}
 \bS-2\nu_{*} \bD= \tau(p_f) \frac{\bS}{|\bS|}.
 \end{equation}
Hence $\bZ:= \bS-2\nu_{*} \bD$ fulfils $|\bZ|=\tau(\pf)$. 
Next, by taking the modulus of \eqref{ZS} it follows 
$$|\bS|-\tau(\pf)=2\nu_*|\bD|.$$
Inserting this back to \eqref{model1}, we get
$$\frac{\bS}{|\bS|}=\frac{\bD}{|\bD|}.$$
Employing this in \eqref{ZS}, we obtain first
$$ \bZ=\bS-2\nu_{*} \bD= \tau(p_f)\frac{\bD}{|\bD|}$$
and then, after taking the scalar product with $\bD$, 
$$ \bZ:\bD=\tau(\pf)|\bD|,$$
which  is the second  assertion in \eqref{model2-2}. 
\par
Next, we assume that  $(\bS, \bD, \pf)$ fulfil \eqref{model2-2}. Then, if  $\bD\neq\bO$,
$$\tau(\pf) |\bD|\leq \bZ:\bD\leq |\bZ||\bD|\leq \tau(\pf)|\bD|, $$ 
which implies 
\begin{equation}\label{ZprodD}
\bZ:\bD=\tau(\pf)|\bD|
\end{equation}
as well as the equality in the Cauchy-Schwarz inequality.
Then necessarilly $$\bZ=a\bD.$$
Inserting this structure in \eqref{ZprodD} we obtain
$$\tau(\pf)|\bD|=a|\bD|^2.$$
Hence $\bZ=\tau(\pf)\frac{\bD}{|\bD|}$ and 
\begin{equation}\label{jos1}
\bS=2\nu_*\bD+\tau(\pf)\frac{\bD}{|\bD|}.
\end{equation}
Also, we have 
$$|\bS|=\left( \frac{2\nu_*|\bD|+\tau(\pf)}{|\bD|}\right)|\bD|=2\nu_*|\bD|+\tau(\pf)$$
which implies (as $\bD\neq \bO$)
$$|\bS|-\tau(\pf)>0, \ \mbox{and also } \frac{\bS}{|\bS|}=\frac{\bD}{|\bD|}.$$
This together with \eqref{jos1} implies  \eqref{model1} for $\bD\neq \bO$.
If $\bD=\bO$ then, by \eqref{model2-2}, $\bS=\bZ$ and $|\bS|\leq \tau(\pf)$ and \eqref{model1} holds. The proof of the equivalence of \eqref{model1} and  \eqref{model2-2} is complete.\par
The proof of the statement (b) is done in the same manner.
 \qed
\\

\section{Definition of weak solution and Main Result}
In order to define the weak solution to the considered problem and to formulate the result, we need to fix the notation. For any $q\in [1,\infty]$ the symbol $\|\cdot\|_q $ stands for the $L^q$-norm in the usual Lebesgue space $L^q(\O)$ while $\|\cdot\|_{1,q}$ for the norm in the usual Sobolev space $W^{1,q}(\O)$. If $X$ is a Banach space of scalar functions then $X^3$ denotes the space of vector-valued functions having three components, each of them belonging to $X$. Similarly $X^{3\times3}$ denotes the space of tensor-valued functions, with each component belonging to $X$. For a Banach space $X$ we denote the relevant Bochner space by $L^q(0,T; X)$. Let us introduce the notation for spaces of solenoidal functions and for  spaces of functions which have zero normal component on the boundary, for the domain $\O$. We set for any $q\in [1,\infty)$
$$\Lnd{q} := \overline{ \left\{ \bv \in C_0^\infty(\Omega)^3; \, \diver \bv =0\right\}}^{\|\cdot\|_q}.$$
Next, we define
\begin{align*}\displ
W_{\mathbf{n}}^{1,2}&  :=  \{\bv\in W^{1,2} (\Omega)^3; \bv\cdot \mathbf{n} = 0 \ \mbox{on}\ \partial \Omega\},\\\displ
W_{\mathbf{n},\diver}^{1,2} & :=  \{\bv\in W^{1,2} (\Omega)^3; \bv\cdot \mathbf{n}=0 \ \mbox{on} \ \partial \Omega; \diver \bv=0\ \mbox{in}\ \Omega\},\\\displ
\Wnm{2} &:= \left ( \Wn{2} \right )^*\!\!, \ \  \  \Wndm{2} := \left ( \Wnd{2} \right )^*\!\!,\\ \displ
W_{\mathbf{n}}^{1,\infty}&  :=  \{\bv\in W^{1,\infty} (\Omega)^3; \bv\cdot \mathbf{n} = 0 \ \mbox{on}\ \partial \Omega\},\\\displ
W_{\mathbf{n},\diver}^{1,\infty} & :=  \{\bv\in W^{1,\infty} (\Omega)^3; \bv\cdot \mathbf{n}=0 \ \mbox{on} \ \partial \Omega; \diver \bv=0\ \mbox{in}\ \Omega\}.
\end{align*}
By the Helmholtz decomposition (in the case that $\O\in C^{1,1}$)  it holds
$$
W^{1,2}_{\bn}=W^{1,2}_{\bn,\diver} \oplus \{\nabla \varphi; \varphi \in W^{2, 2}(\Omega), \nabla \varphi \cdot \bn =0 \textrm { on } \partial \Omega\}.
$$
Note that such a decomposition is not valid for $(W^{1,2}_{0}(\Omega))^3$.  \par
The symbol $\bD\boldsymbol\varphi$ stands for the symmetric part of the gradient of a vector-valued function $\boldsymbol\varphi$, i.e. $\bD\boldsymbol\varphi:=\frac{\nabla\boldsymbol\varphi+(\nabla\boldsymbol\varphi)^T}{2}$.\par
In what follows, we also set for simplicity and without loss of any generality
$$\varrho_* = 2\nu_* = \gamma_* = K=q_*=1.$$ 
\begin{defn}[Definition of weak solution]
Let $s_*>0$, 
\be\label{def-ws}
\begin{split}
&\bv_0\in \Lnd{2}, \ p{_0}\in L^\infty(\O), \ \bb\in L^2(0,T; \Wnm{2}),\\
\end{split}\ee
and one of the following requirements be satisfied 
\begin{alignat}{1}\label{1}
& { p_s\in L^\infty(Q_T), p_s(0) \in L^\infty(\O),} g\in L^q(Q_T),  \partial_t p_s-\Delta p_s\in L^q(Q_T) \ \mbox{ with } q>\frac{5}{2},\\
& p_s\in L^q(0,T; W^{1,q}(\O)) \mbox{ with  } q>10 \ \mbox{and } g\in L^q(Q_T) \ \mbox{ with } q>\frac{5}{2}.\label{2}
\end{alignat}
 We say that $(\bv,\pf, p, \bS, \bs)$ is  a weak solution to the problem \eqref{system}-\eqref{g2} if
\begin{align*}
&\bv \in  L^\infty(0,T; \Lnd{2})\cap L^2(0,T;W^{1,2}_{{\bf n},{\rm div}}),  \ \partial_t\bv\in (L^2(0,T;W^{1,2}_{\bn}(\O))\cap L^{5}(Q_T)^3)^*,\\
&\pf \in L^\infty(Q_T)\cap L^2(0,T; W^{1,2}(\O)), \ \partial_t \pf \in (L^2(0,T; W^{1,2}(\O)))^*,\\
& p=p_1+p_2 \mbox{ where } p_1\in L^2(Q_T) \mbox{ and } p_2 \in L^{\frac{5}{4}}(0,T; W^{1, \frac{5}{4}}(\O)), \\
&\bS\in  L^{2}(Q_T)^{3\times 3}, \ \bs\in L^2(\Sigma_T)^3,  \\
&\langle \partial_t\bv, \bw\rangle + (\bS, \bD\bw)
 { + (\diver(\bv\otimes \bv) , \bw)} + (\bs, \bw_{\btau})_{\partial \Omega}
 = \langle\bb,  \bw\rangle + (p_1,\diver \bw)-(\nabla p_2, \bw)\\
 &\qquad \textrm{ for all } \bw\in W^{1,2}_{\bn}  \textrm{ and a.e. in} \ (0,T), \\
& \langle \partial_t \pf, z\rangle - ( \pf \bv, \nabla z )+ (\nabla \pf,\nabla z)=(g, z)-(p_s\bv, \nabla z)\\
&\qquad \mbox{for all} \ z\in W^{1,2}\  \mbox{and a.e. in } (0,T),\\
& \bD\bv =  \frac{\left(|\bS| - \tau(\pf)\right)^{+}}{|\bS|} \bS \ \ \mbox{where}\ \ \tau({\pf}) = (p_s- \pf)^+ \  \textrm{ for a.a. }
(t,x) \in Q_T, \\
&\bv_\tau=\frac{(|\bs|-s_*)^+}{|\bs|}\bs \  \textrm{ for a.a. }
(t,x) \in \Sigma_T,\\
&\lim_{t\rightarrow 0+}\|\bv(t)-\bv_0\|_2=0 \  \mbox{and }\  \lim_{t\rightarrow 0+}\|\pf(t)-p_0\|_2=0.
\end{align*}
\end{defn}

\begin{thm}[Main Theorem]\label{main-thm}
For any $\O\in C^{1,1}$, $T > 0$ and for arbitrary  $ \bv_0, p{_0}, p_s, \bb$ fulfilling \eqref{def-ws} and for arbitrary $g$ and $p_s$ fulfilling either \eqref{1} or \eqref{2}, there exists a weak solution to the problem \eqref{system} in the sense of Definition \ref{def-ws}. 
\end{thm}
\begin{rem}
We wish to emphasize that due to Proposition \ref{equivalent-formulation}, the tensorial constitutive equation $\bD\bv=\frac{(|\bS|-\tau(\pf))^+}{|\bS|}\bS$ in $Q_T$ as well as the vectorial equation $\bv_\tau=\frac{(|\bs|-s_*)^+}{|\bs|}\bs$ on $\Sigma_T$ can be replaced by any of its equivalent forms. It is in particular interesting that the tensorial equations can be characterized by two (scalar) inequalities.
\end{rem}
Note that Theorem \ref{main-thm} presents the existence result to a supercritical problem; this is a problem where the solution itself is not an admissible test function in the weak formulation of the governing equations. Indeed, in our case $\bv$ belongs to $L^{\frac{10}{3}}(Q_T)^3$, however, admissible test functions have to be from $L^5(Q_T)^3$ due to the fact that $\diver(\bv\otimes\bv)=\sum_{r=1}^3 v_r\frac{\partial \bv}{\partial x_r}$ and $\nabla p_2$ belong to $L^{\frac{5}{4}}(Q_T)^3$. This is the reason why we cannot involve in our analysis such tools as the energy equality, used in the analysis of planar time-dependent flows in Chupin and Math\'e \cite{CM} or the higher differentiability techniques used in \cite{LS2000} and \cite{She}, also in the analysis of two-dimensional unsteady flows of the Bingham fluids and in the analysis of steady flows in three dimensions in \cite{MaRuSh}. Neither can we incorporate the tools of calculus of variations suitable for Stokes-type problems (see for example \cite{FS} and the references therein). On the other hand, we intentionally aim at avoiding tools such as multivalued calculus  or variational inequalities \cite{Aubin} in our analysis, see \cite{DL}, or  \cite{Rod} (which is considered however in a different context). In our opinion, the concept of solution (expressed in terms of identities) considered here is stronger, its large-data existence can be proved and has some other advantages. For example, it forms the foundation for a direct application of mixed finite element (or spectral) methods.
In order to identify the non-linear constitutive equation pointwise in the considered domain $(0,T)\times\O$ when taking the limit from the approximative problem to the original one, and in order to overcome difficulties connected with the low integrability of $\bv$, we incorporate the so-called $L^\infty$-truncation method. This method replaces $\bv^n-\bv$, where $\{\bv^n\}_{n=1}^{+\infty}$ is solution of a suitably constructed approximative problem, by a truncated function that coincides with $\bv^n-\bv$ on a large set and the measure of the complementary set can be made arbitrarly small uniformly with respect to $n$. Although the origin of the method goes back to elliptic problems with an $L^1$-right-hand side (see \cite{DM1998}, \cite{FMS1997} and \cite{R1997}), we refer here mainly to its development for evolutionary problems in fluid mechanics, see \cite{FMS}, \cite{BMR}, \cite{Wolf}.
The result by Wolf \cite{Wolf} similarly as those by Solonnikov (see \cite{Sol1} and \cite{Sol2}) and Koch and Solonnikov \cite{KochSol} concerning the properties of evolutionary Stokes-like systems with no-slip boundary conditions indicate the difficulties connected with the impossibility to establish the integrability of the pressure $p$ for generalizations of the Navier-Stokes equations (with variable viscosity) in three-dimensions. This is why we treat the stick-slip boundary conditions in this study. It reveals that the analysis of the three-dimensional evolutionary supercritical problems associated with the stick-slip boundary conditions differs remarkably from the analysis of analogous problems connected with the no-slip boundary conditions. We refer to \cite{BM} for a detailed discussion of this issue noting that Theorem \ref{main-thm} guarantees that $p\in L^1(Q_T)$. We remark that the integrability of the pressure is important in the analysis of problems with the viscosity dependent on the temperature (see \cite{BFM2009}, \cite{BM17} or \cite{MZ}) or the viscosity dependent on the pressure (see \cite{BMR} or \cite{BMR2009}), but it is also an interesting mathematical question itself.

\section{Approximations}
Before introducing the approximations, we recall that in Section 3 we set
$$ \bZ=\bS-\bD\bv,$$
and analogously we can also define 
$$\bz:=\bs-\bv_\tau.$$
 For any $n\in \N$, let $G_n:\R\to \R$ be a smooth function such that $G_n(u)=1$ if $|u|\leq n$, $G_n(u)=0$ if $|u|\geq 2n$ and $|G'_n|\leq \frac{2}{n}$. We consider the following approximative system:
\be\label{system22}
\ba{l}\displ
\diver\bv = 0 \textrm{ in } Q_T, \\ \displ
\pa_t\bv + \diver(\bv\otimes \bv) G_n(|\bv|^2)-  \diver\bD \bv -\diver  \bZ +\nabla p =  \bb \textrm{ in } Q_T, \\\displ
\pa_t \pf+ \bv \cdot \na \pf-  \Delta \pf=g+\bv\cdot\nabla p_s \textrm{ in } Q_T, \\\displ
\bZ=\mathcal{Z}_n(\pf, \bD\bv) := (p_s-\pf)^+  \frac{\bD\bv}{|\bD\bv|+\frac{1}{n}}\textrm{ in } Q_T, \\\displ
\bz= \zeta_n(\bv_\tau) := s_* \frac{\bv_{\tau}}{|\bv_{\tau}|+\frac{1}{n}} \  \ \mbox{on} \ \ \Sigma_T,\\\displ
\bv\cdot\bn=0 \ \mbox{and} \ \na \pf \cdot \bn=0 \ \mbox{on} \ \ \Sigma_T,\\\displ
\bv(0)=\bv_0 \ \mbox{and} \ \pf(0)={p}_0 \ \mbox{in} \ \O.
\ea
\ee
It is not difficult to check that if \ $\bZ=\mathcal{Z}_n(\pf, \bD)$ and $\hat{\bZ}=\mathcal{Z}_n(\pf, \hat{\bD})$, then 
\be\label{Prop-monotonicity} (\bZ-\hat{\bZ}):(\bD-\hat{\bD})\geq \frac{(p_s-\pf)^+ }{n} \frac{(|\bD|-|\hat{\bD}|)^2}{\left(|\bD|+\frac{1}{n}\right)\left(|\hat{\bD}|+\frac{1}{n}\right)}\geq 0.
\ee
A similar monotone property holds for $\bz= \zeta_n(\bv_\tau)$.
\begin{prop}\label{prop}
Let $n\in \N$ be fixed and $s_*>0$. Let $ \bv_0\in \Lnd{2}$, $p_0\in L^2(\O)$, $\bb\in L^2(0,T; \Wnm{2})$,  $g\in L^2(Q_T)$  and $p_s\in L^{{ 5}}(Q_T)$, then there exists a weak solution to the problem \eqref{system22}, i.e.  a quadruple $(\bv^n,\pf^n, \bZ^n, \bz^n)$ such that
\begin{align}
&\bv^n \in  L^\infty (0,T; \Lnd{2})\cap L^2(0,T;W^{1,2}_{{\bf n},{\rm div}}), \ \ \partial_t\bv^n\in (L^2(0,T;W^{1,2}_{\bn}(\O))^*, \label{an1}\\
&\pf^n \in  L^\infty (0,T; L^2(\O))\cap L^2(0,T; W^{1,2}(\O)),\ \ \partial_t \pf^n \in (L^4(0,T; W^{1,2}(\O)))^*,\label{an2}\\
&\bZ^n\in  L^{\frac{10}{3}}(Q_T)^{3\times 3}, \ \ \bz^n\in L^\infty(\Sigma_T)^3, \label{an3}\\
&\label{an4}\ba{l}\displaystyle\vspace{4pt}
\int_0^T\!\!\! \langle \partial_t\bv^n, \bw\rangle dt + \int_{Q_T}\!\!\!\! \bD \bv^n \!:\! \bD \bw + \bZ^n\!:\! \bD\bw
 + G_n(|\bv^n|^2)\diver(\bv^n\otimes \bv^n)\! \cdot \!\bw \,dx dt \\ \displaystyle\vspace{4pt}
 + \int_{\Sigma_T}\!\!\!\!\! \bz^n\!\cdot\! \bw_{\btau}+ \bv_\tau^n\!\cdot\! \bw_{\btau} \,d\sigma_xdt
  = \int_0^T\!\!\! \langle\bb,  \bw\rangle dt\ \  \textrm{ for all } \bw\in L^2(0, T; W^{1,2}_{\bn, \diver}), 
 \ea 
 \\&\label{an5} \ba{l}\displaystyle\vspace{4pt}
\int_0^T\!\!\!   \langle \partial_t \pf^n, z\rangle \,dt { - } \int_{Q_T}\!\!\! \pf^n \bv^n\!\cdot\! \nabla z  + \nabla \pf^n\!\cdot\!\nabla z\,dxdt= \int_{Q_T}\!\!\! \!g z- p_s\bv^n\!\cdot\!\nabla z  \, dxdt \\ \displaystyle\vspace{4pt}
\hfill\mbox{ for all } z\in L^4(0,T; W^{1,2}(\O)),
  \ea \\
& \bZ^n=\mathcal{Z}_n(\pf^n, \bD\bv^n) \textrm{ a.e. in } Q_T, \label{an6}\\ 
&\bz^n= \zeta_n(\bv_\tau^n) \textrm{  a.e. in } \Sigma_T, \label{an7}\\
&\lim_{t\rightarrow 0+}\|\bv^n(t)-\bv_0\|_2=0 \ \mbox{and }\  \lim_{t\rightarrow 0+}\|\pf^n(t)-p{_0}\|_2=0. \label{an8}
\end{align}
\end{prop}

\Pr
Due to the presence of $G_n$ that truncates the convective term and the properties of the approximations $\mathcal{Z}_n$ and $ \zeta_n$ introduced above, the proof of the existence of weak solutions to the problem \eqref{system22} is a variant of the standard monotone operator technique (see \cite{Lions}, \cite{Lady} or \cite{BMR}). To be more specific, we briefly outline the proof using the Galerkin method. 
Since $n$ is fixed, we write $(\bv,\pf, \bZ, \bz)$ instead of $(\bv^n,\pf^n, \bZ^n, \bz^n)$ in the proof.
\par
{\textit{Step 1. Galerkin system.}
 Let $\{\mathbf{w}^i\}_{i\in\N}$ be an orthogonal basis in $W_{\mathbf{n}, \diver}^{1,2}$ consisting of eigenfunctions of the Stokes operator subject to $\bv\!\cdot\! \bn\!=\!0$ and $[(\bD\bv) \bn]_{\tau}\!=\!\bnul$ on $\Sigma_T$. Let analogously $\{z_j\}_{j\in \mathbb{N}}$ be an orthogonal basis in $W^{1,2}(\Omega)$ consisting of  eigenfunctions of the Laplace operator subject to the relevant homogeneous boundary conditions. Then the local in time existence of 
 \be
 \bv^{m}(t,\bx):=\sum_{r=1}^m c_r^{m}(t) \mathbf{w}^r(\bx), \ \ \ p_{f}^{m}(t,\bx):=\sum_{r=1}^m d_r^{m}(t) z^r(\bx)
\ee
satisfying 
\be\label{alpha}\ba{l}\displ
\left(\! \frac{d \bv^{m}}{dt}, \bw^r\!\!\right)\!+\! (\bD \bv^{m}\!, \bD \bw^r) 
\!+\! ( \mathcal{Z}_n(\pf^m, \bD\bv^m\!), \bD \bw^r ) \!+ \!( \diver (\bv^{m}\!\otimes\!\bv^{m}\!)G(|\bv^{m}|^2),\bw^r)
\\\displ \hfill+(\bv^m_\tau, \bw^r)_{\pa\O}+ ( \zeta_n(\bv_\tau^{m}), \bw^r)_{\pa\O} = \langle\bb,\bw^r\rangle, \ \ \ r=1,\dots, m,
\ea\ee
and 
\be\label{beta}\ba{l}\displ
\!\!{ (\partial_t p_{f}^{m}, z^r )} \!-\!  ( p_{f}^{m}\bv^{m},\nabla z^r ) \!+\!  ( \nabla p_{f}^{m},  \nabla z^r) = (g, z^r)\!-\!(p_s\bv^{m},\nabla z^r ), \ \ \ r=1,\dots, m,
\ea\ee
together with the corresponding initial conditions $\bv^m_0$ and $p_0^m$, obtained by projecting $ \bv_0\in \Lnd{2}$ onto the span of $[ \mathbf{w}^1,\dots ,\mathbf{w}^m]$ and  $p_0\in L^2(\O)$ onto the span of $ [z^1,\dots, z^m]$, follows from the Caratheodory theory for systems of ordinary differential equations. 
\par Global in time existence is, as usual, a consequence of the uniform estimates which we show next. 
\par \textit{Step 2. Uniform estimates.} Multiplying \eqref{alpha} by $c^m_r(t)$ and \eqref{beta} by $d^m_r(t)$ and taking the sum over $r$ from $1$ to $m$, we obtain
\begin{alignat}{1}
\label{gamma}
&\ba{l}\displ
\frac{1}{2} \frac{d}{dt}\|\bv^{m}(t)\|_2^2 +   {\|\bD \bv^{m}(t)\|}_2^2
+ \int_{\Omega}\!\! { \mathcal{Z}_n(\pf^m, \bD\bv^m\!)}\!:\!\bD\bv^m\,dx 
+\|\bv^m_\tau(t)\|_{2,\partial \Omega}^2  \\ \displ 
\hfill +\int_{\partial \Omega}\!\!\! \zeta_n(\bv_\tau^{m})\! \cdot\! \bv_{ \tau}^{m}\,d\sigma=\langle \bb, \bv^{m}\rangle,
\ea
\\
\label{delta}
&\frac{1}{2} \frac{d}{d t} \|\pf^{m}(t)\|_2^2 +  \| \nabla \pf^{m}(t)\|_2^2 = (g, \pf^m)-(p_s\bv^m, \nabla \pf^m), 
\end{alignat}
By Korn's and Young's inequalities (see for example \cite[Lemma 1.11]{BMR} for details), using also the fact that the last two terms at the right-hand side of \eqref{gamma} are non-negative, one concludes from \eqref{gamma} that 
\begin{equation}\label{AEO}\ba{l}\displ
\sup_{t\in [0,T]} \|\bv^{m} (t)\|_2^2 +  \int_{Q_T}\!\!\!\! |\bD \bv^{m}|^2 { + |\bv^{m}|^{10/3}\,dxdt} + \int_{\Sigma_T} \!\!\! |\bv^m_\tau|^2 d\sigma dt\\ \displ
\hfill \leq  C\|\bb\|_{L^2(0,T;  \Wnm{2})}^2+\|\bv_0\|_2^2=:C(\bb, \bv_0),
\ea\end{equation}
{ where we also used the interpolation inequality
\begin{align}
&\|z\|_{\frac{10}{3}}\leq \|z\|_2^{\frac{2}{5}}\|z\|_6^{\frac{3}{5}}\leq C\|z\|_2^{\frac{2}{5}}\|z\|_{1,2}^{\frac{3}{5}}.\label{dt0}
\end{align}}
Similarly, using also { 
$$
\int_0^T \!\!\!\! |(p_s\bv^m, \nabla \pf^m)| \, dt \le \|p_s\|_{L^{5}(Q_T)} \|\bv^m\|_{L^{10/3}(Q_T)} \|\nabla \pf^m\|_{L^2(Q_T)}, 
$$}
one obtains from \eqref{delta}, using also \eqref{AEO}, that
\begin{equation}\label{AE1}
\sup_{t\in [0,T]}\|\pf^{m}(t)\|_2^2+ \int_{Q_T}\!\!\!\! |\nabla \pf^{m}|^2\,dxdt \leq C\|g\|_{L^2(Q_T)}+ C(\bb,\bv_0)\|p_s\|_{L^{{5}}(Q_T)}+ \|{p}_0 \|_2^2.
\end{equation}
By the interpolation inequalities \eqref{dt0} and 
\begin{align}
& \|z\|_{4}\leq \|z\|_2^{\frac{1}{4}}\|z\|_6^{\frac{3}{4}}\leq C\|z\|_2^{\frac{1}{4}}\|z\|_{1,2}^{\frac{3}{4}},\label{inter*}
\end{align}
and by the trace inequalities (see \cite[Lemma { 1.11}]{BMR}), we obtain 
\begin{equation}\label{AE2}
\sup_m \left(\!\|\pf^m\|_{\frac{10}{3}, Q_T}+\|\bv^m_\tau\|_{\frac{8}{3}, \Sigma_T}\!\!\right)<+\infty,
\end{equation}
and also 
\begin{equation}\label{AE2a}
\sup_m \left(\!\!\int_0^T\!\!\! \left(\|\bv^m(t)\|_4^{\frac{8}{3}} + \|\pf^m\|_4^{\frac{8}{3}}\right)dt\!\!\right)<+\infty.
\end{equation}
It then follows from the explicit formulas for $\mathcal{Z}_n$ and $\zeta_n$ that $\bZ^m:=\mathcal{Z}_n(\pf^m, \bD\bv^m\!)$ and $\bz^m:=\zeta_n(\bv_\tau^m)$ fulfil 
\begin{equation}\label{AE4}
\sup_m\left(\|\bZ^m\|_{\frac{10}{3}, Q_T} +\|\bz^m\|_{\infty, \Sigma_T}\right)<+\infty.
\end{equation}
Finally, the fact that the projectors
$$W^{1,2}_{\bn, \diver}\longmapsto [ \mathbf{w}^1,\dots ,\mathbf{w}^m],\ W^{1,2}(\O) \longmapsto  [z^1,\dots, z^m]
$$
are continuous and \eqref{AE2a} imply that
\begin{equation}\label{AE5}
\sup_m \!\left(\!\|\partial_t\bv^m\|_{L^2(0,T; \Wnm{2})}+\|\partial_t\pf^m\|_{L^{\frac{4}{3}}(0,T; W^{-1,2})}\!\right)\!<+\infty.
\end{equation}
\par \textit{Step 3. Limit.} The above uniform estimates imply the existence of $\bv$, $\pf$, $\bZ$ and $\bz$ and subsequences of $\{\bv^m\}, \{\pf^m\}, \{\bZ^m\}$ and $\{\!\bz^m\!\}$ converging weakly (or *-weakly) to $\bv, \pf, \bZ$ and $\bz$ in the function spaces indicated in Proposition \ref{prop}, and fulfilling the following strong convergences (due to Aubin-Lions compactness lemma and its variant, see \cite[Lemma 1.12]{BMR}, involving the trace theorem):
\begin{align}
&\bv^m\to \bv \mbox{ a.e. in } Q_T \mbox{ and strongly in } L^q(Q_T)^3 \mbox{ for any } q\in \left[1, \frac{10}{3}\right),\label{L1}\\
&\pf^m\to \pf  \mbox{ a.e. in } Q_T \mbox{ and strongly in } L^q(Q_T)^3 \mbox{ for any } q\in \left[1, \frac{10}{3}\right),\label{L2}\\
& \bv^m_\tau \to \bv_\tau  \mbox{ a.e. in } \Sigma_T \mbox{ and strongly in } L^q(\Sigma_T)^3 \mbox{ for any } q\in \left[1, \frac{8}{3}\right).\label{L3}
\end{align}
These weak and strong convergences suffice to show that $\bv, \pf, \bZ$ and $\bz$ fulfil the weak formulations {} \eqref{an4}--\eqref{an5}} stated in Proposition \ref{prop}.
\par Since the proof of the attainment of the initial conditions is standard, see e.g. \cite{MR2005}, it remains to show that $\bZ=\mathcal{Z}_n(\pf, \bD\bv)$ and $\bz=\zeta_n(\bv_\tau)$.
\par \textit{Step 4. Attainment of the constitutive equations.} We first notice that \eqref{L3} together with \eqref{AE4} imply, by Lebesgue's theorem that
$$ \bz^m=\zeta_n(\bv_\tau^m) \rightharpoonup \zeta_n(\bv_\tau) \mbox{ weakly in } L^2(\Sigma_T)^3.$$
It implies that 
\begin{equation}\label{L4}
\bz=\zeta_n(\bv_\tau) \mbox{ a.e. in } \Sigma_T
\end{equation}
and
\begin{equation}\label{L5}
\lim_{m\to+\infty} \int_{\Sigma_T}\! \!\!\bz^m\!\cdot\!\bv^m_\tau d\sigma dt= \int_{\Sigma_T}\! \!\! \zeta_n(\bv_\tau)\!\cdot\! \bv_\tau d\sigma dt.\end{equation}
Next, integrating \eqref{gamma} over $(0,T)$ and taking limsup of the resulting identity, we obtain, using the above convergences and the weak lower semicontinuity of the $L^2$-norm, that 
\begin{equation}\label{EI1}\begin{split}
&\frac{1}{2} \|\bv(t)\|_2^2+\int_{Q_T}\!\!\!\! |\bD\bv|^2dxdt+\int_{\Sigma_T}\!\!\!\! |\bv_\tau|^2d\sigma dt+ \int_{\Sigma_T}\!\!\!\!\zeta_n(\bv_\tau)\!\cdot\!\bv_\tau d\sigma dt\\
&+\limsup_{m\to+\infty}\int_{Q_T}\!\!\!\!\!\bZ^m\!\!:\!\bD\bv^mdxdt\leq \int_0^T\!\!\!\langle \bb, \bv\rangle dt +\frac{1}{2}\|\bv_0\|_2^2.
\end{split}\end{equation}
On the other hand, taking $\bw=\bv$ in the established weak formulation of the equation for $\bv$, we get, using also \eqref{L4}, 
\begin{equation}\label{EI2}\begin{split}
&\frac{1}{2} \|\bv(t)\|_2^2+\int_{Q_T}\!\!\!\! |\bD\bv|^2dxdt+\int_{\Sigma_T}\!\!\! |\bv_\tau|^2d\sigma dt+ \int_{\Sigma_T}\!\!\!\!\zeta_n(\bv_\tau)\!\cdot\!\bv_\tau \,d\sigma dt\\
&+\int_{Q_T}\!\!\!\!\! \bZ\!:\!\bD\bv \,dxdt = \int_0^T\!\!\!\langle \bb, \bv\rangle dt +\frac{1}{2}\|\bv_0\|_2^2.
\end{split}\end{equation}
Comparing \eqref{EI1} with \eqref{EI2}, we conclude that
\begin{equation}\label{EI3}
\limsup_{m\to+\infty} \int_{Q_T}\!\!\!\!\!\bZ^m\!\!:\!\bD\bv^mdxdt\leq \int_{Q_T}\!\!\!\!\! \bZ\!:\!\bD\bv\,dxdt.
\end{equation}
Finally, it follows from \eqref{Prop-monotonicity} that
\begin{equation}\label{M1}
0\leq \!\int_{Q_T}\!\!\!\!\!\!\left(\!\mathcal{Z}_n(\pf^m, \bD\bv^m\!)- \mathcal{Z}_n(\pf^m, \bA\!)\right)\!:\!\left(\bD\bv^m-\bA\right)dxdt\ \ \mbox{   for all } \bA \in L^2(Q_T).
\end{equation}
Since, by \eqref{L2},
$$  \mathcal{Z}_n(\pf^m, \bA\!):=(p_s - \pf^m)^+\!\frac{\bA}{|\bA|+\frac{1}{n}}\to(p_s-pf)^+\!\frac{\bA}{|\bA|+\frac{1}{n}}=:\mathcal{Z}_n(\pf, \bA\!) \mbox{ strongly in } L^2(Q_T) $$
and 
$$\bD\bv^m\rightharpoonup\bD\bv \ \mbox{ weakly in } L^2(Q_T), $$
we conclude from \eqref{M1} and \eqref{EI3} that 
\begin{equation}\label{M2}
0\leq \!\int_{Q_T}\!\!\!\!\!\!\left(\bZ- \mathcal{Z}_n(\pf, \bA\!)\right)\!:\!\left(\bD\bv-\bA\right)dxdt\ \ \mbox{   for all } \bA \in L^2(Q_T).
\end{equation}
The choice $\bA=\bD\bv\pm\lambda\bB$ for $\bB\in L^2(Q_T)$ arbitrary and $\lambda>0$, leads to
$$  0\leq \pm \!\int_{Q_T}\!\!\!\!\!\!\left(\bZ- \mathcal{Z}_n(\pf,\bD\bv\pm\lambda\bB\!)\right)\!:\!\bB\, dxdt\ \ \mbox{   for all } \bB \in L^2(Q_T).$$
Letting $\lambda\to 0^+$, we obtain
$$ 0= \!\int_{Q_T}\!\!\!\!\!\!\left(\bZ- \mathcal{Z}_n(\pf,\bD\bv\!)\right)\!:\!\bB\, dxdt\ \ \mbox{   for all } \bB \in L^2(Q_T),$$
which implies $\bZ=\mathcal{Z}_n(\pf,\bD\bv\!)$ a.e. in $Q_T$.\\ The proof of Proposition \ref{prop} is complete.
\qed

\begin{prop}\label{prop52}
Let all the assumptions in Proposition \ref{prop} be satisfied. In addition, assume that { $p_0\in L^{\infty}(\Omega)$} and one of the following requirements holds:
\begin{alignat}{1}\label{ass1}
& p_s(0) \in L^\infty(\O){, \ p_s \in L^{\infty}(Q_T)} \ \mbox{ and } g, \partial_tp_s-\Delta p_s\in L^q(Q_T) \ \mbox{ with } q>\frac{5}{2},\\
&p_s\in L^q(0,T; W^{1,q}(\O)) \mbox{ with  } q>10 \ \mbox{and } g\in L^q(Q_T) \ \mbox{ with } q>\frac{5}{2},\label{ass2}
\end{alignat}
then, for each $n\in\N$,  there exists a weak solution to the problem \eqref{system22} in the sense of  Proposition \ref{prop} satisfying $\pf^n\in L^\infty(Q_T)$. In fact,
$$\sup_n \|\pf^n\|_{L^\infty(Q_T)}<+\infty.$$
Consequently, $\partial_t\pf^n\in L^2(0, T; W^{-1,2}(\O))$ and \eqref{an5} holds for all $z\in L^2(0,T;W^{1,2}(\O))$.
\end{prop}
\Pr
In what follows we shall prove explicitly that $\pf\in L^\infty (Q_T)$ using the Moser iteration technique.
By the interpolation inequality  \eqref{dt0},
it follows that
\begin{equation}\label{dt1}
\int_0^T \|z\|_{\frac{10}{3}}^{\frac{10}{3}}\leq C\left( \!\sup_{t\in [0,T]}\|z\|_2\!\!\right)^{\frac{4}{3}}\int_0^T \|z\|_{1,2}^{2}.
\end{equation}
Consequently,
\begin{equation}\label{dt}
\pf\in L^{\frac{10}{3}}(Q_T) \ \ \mbox{and} \ \ \bv\in L^{\frac{10}{3}}(Q_T)^3.
\end{equation}
Let us first consider the case given by \eqref{ass1}. Then, once we set $h:=-(\partial_tp_s-\Delta p_s)$, $G:=g+h$ and $P:=\pf-p_s$, we can rewrite the third equation in \eqref{system22} as
\begin{equation}\label{E1}
\partial_tP+\bv\cdot\nabla P-\Delta P=G \ \ \mbox{with} \ G\in L^q(Q_T) \ \mbox{and} \ q>\frac{5}{2}.
\end{equation}
 For $s>2$ and $m\in\N$ consider  $|P_m|^{s-2}P_m$ with $P_m:=T_m(P)$ as test function in the weak formulation of \eqref{E1}. Here $T_m:\R\to\R$ is defined through $T_m(z)=z$ if $|z|\leq m$ and $T_m(z)=m {\ \mathrm{sgn}} \, z $ if $|z|>m$. Note that $|P_m|^{s-2}P_m$ is an admissible test function. After integrating by parts and employing  $\diver \bv =0$, we get 
\begin{equation}
\frac{1}{s}\frac{d}{dt}\|P_m\|_s^s+(s-1)\!\!\int_\O\! |\nabla P_m|^2|P_m|^{s-2}\leq \int_\O |G| |P_m|^{s-1}.
\end{equation}
Next, integrating with respect to the time, straightforward computations imply
\begin{equation}\label{pfs}
\| |P_m(t)|^{\frac{s}{2}}\|_2^2+\frac{4s(s-1)}{s^2} \| \nabla |P_m|^{\frac{s}{2}}\|_{2,Q_T}^2\leq s\! \int_{Q_T} |G| |P_m|^{s-1}
+ \|P(0)\|_s^s=:A.
\end{equation}
Since $\frac{4s(s-1)}{s^2}>1$, it follows from \eqref{pfs} that
\begin{equation}\label{A}
\sup_{t\in [0,T]} \| |P_m(t)|^{\frac{s}{2}}\|_2\leq A^{\frac{1}{2}}, \ \ \int_0^T \| \nabla |P_m(t)|^{\frac{s}{2}}\|_2^2\,dt\leq A.
\end{equation}
Using \eqref{dt0} and \eqref{dt1} with $z=|P_m|^{\frac{s}{2}}$, and combining the result with \eqref{A}, we obtain
\begin{equation}\label{53s}\begin{split}
\int_0^T \|P_m(t)\|_{\frac{5s}{3}}^{\frac{5s}{3}}dt &= \int_0^T \||P_m(t)|^{\frac{s}{2}}\|_{\frac{10}{3}}^{\frac{10}{3}}dt \\ 
&\leq C \left(\!\sup_{t\in [0,T]} \!\!\| |P_m(t)|^{\frac{s}{2}}\|_2 \right)^{\frac{4}{3}} \int_0^T \| \nabla |P_m(t)|^{\frac{s}{2}}\|_2^2 dt \leq C  A^{\frac{5}{3}}.
\end{split}
\end{equation}
The definition of $A$ and \eqref{53s} then leads to
\begin{equation}
\| P_m\|_{\frac{5s}{3}, Q_T}\leq s^{\frac{1}{s}} C^{\frac{3}{5s}}\| G\|_{q, Q_T}^{\frac{1}{s}}\| P_m\|_{q'(s-1), Q_T}^{\frac{s-1}{s}}+C^{\frac{3}{5s}}\|P(0)\|_\infty.
\end{equation}
We can introduce the following iteration scheme. Setting  
\begin{equation}\label{rel1}
s_0:=\frac{10}{3}, \ \frac{q}{q-1}(\tilde{s}_i-1):=s_i \ \mbox{ and } s_{i+1}:= \frac{5}{3}\tilde{s}_i,
\end{equation}
which leads to $\tilde{s}_{i}= \frac{q-1}{q} s_i+1$ and hence
\begin{equation}\label{rel2}
s_{i+1}=\frac{5}{3}\frac{q-1}{q}s_i+\frac{5}{3}, \ \tilde{s}_{i+1}=\frac{5}{3}\frac{q-1}{q} \tilde{s}_i+1,
\end{equation}
we obtain
$$ \|P_m\|_{s_{i+1}, Q_T}\leq \tilde{s}_i^{\frac{1}{\tilde{s}_i}} C^{\frac{3}{5\tilde{s}_i}}\| G\|_{q, Q_T}^{\frac{1}{\tilde{s}_i}}\|P_m\|_{s_i, Q_T}^{\frac{\tilde{s}_i-1}{\tilde{s}_i}}+C^{\frac{3}{5\tilde{s}_i}}\|P(0)\|_\infty. $$
Noticing  that  $$\frac{5}{3}\frac{q-1}{q}>1 \iff q>\frac{5}{2},$$
we observe that $s_i\to+\infty$ as $i\to+\infty$. By iteration, we get
\be\label{product}\ba{c}\displ
\!\!\!\!\!\!\|P_m\|_{s_{i+1}, Q_T}\!\\ 	\displ\leq\! C^{\sum_{j=0}^i\!\frac{3}{5\tilde{s}_j}\!\prod_{h=j+1}^{i} \!\!\frac{\tilde{s}_h-1}{\tilde{s}_h}}\!\!\prod_{j=0}^i \tilde{s}_j^{ \frac{1}{\tilde{s}_j}\prod_{h=j+1}^{i}\frac{\tilde{s}_h-1}{\tilde{s}_h}}\!
\|G\|_q^{\sum_{j=0}^i \frac{1}{\tilde{s}_j}\! \prod_{h=j+1}^{i}\frac{\tilde{s}_h-1}{\tilde{s}_h}\! } \|P\|_{s_0} ^{\prod_{j=0}^i \frac{\tilde{s}_j-1}{\tilde{s}_j}}\\ \displ
\!\!\!\!\!\!\!\!+\!\sum_{j=1}^i C^{\sum_{r=j}^i\frac{3}{5\tilde{s}_r}\prod_{h=r+1}^{i} \!\!\frac{\tilde{s}_h-1}{\tilde{s}_h} }\prod_{k=j+1}^i \tilde{s}_k^{ \frac{1}{\tilde{s}_k}\!\prod_{h=k+1}^{i}\!\frac{\tilde{s}_h-1}{\tilde{s}_h}}\\ \displ\cdot
\|G\|_q^{\sum_{k=j+1}^i \frac{1}{\tilde{s}_k} \prod_{h=k+1}^{i}\frac{\tilde{s}_h-1}{\tilde{s}_h} }
\!\|P(0)\|_\infty^{\prod_{k=j+1}^i \!\frac{\tilde{s}_k-1}{\tilde{s}_k}}
\ea
\ee
Next, we use
\begin{align*}
& \|G\|_{q, Q_T}\leq C_1 \ \mbox{ where } C_1:=\max \{ 1, \|G\|_{q, Q_T}\},\\
&\|P(0)\|_\infty\leq C_2 \ \mbox{ where } C_2:=\max \{ 1, \|P(0)\|_{\infty}\}.
\end{align*}
Since $\frac{\tilde{s}_h-1}{\tilde{s}_h}\leq 1$, we notice that all products can be bounded by $1$. Consequently, \eqref{product} leads to (assuming that $C\geq 1$)
\be\ba{l}\displ
\|P_m\|_{s_{i+1}, Q_T}\leq C^{\sum_{j=0}^i\frac{3}{5\tilde{s}_j}} e^{\sum_{j=0}^i \frac{\log\tilde{s}_j}{\tilde{s}_j}} C_1^{\sum_{j=0}^i \frac{1}{\tilde{s}_j}}\max\{1, \|P_m\|_{s_0}\}\\
\displ\hfill +\sum_{j=0}^i C^{\sum_{r=j}^i\frac{3}{5\tilde{s}_r}} e^{\sum_{r=j+1}^i \frac{\log\tilde{s}_r}{\tilde{s}_r}} C_1^{\sum_{r=j+1}^i \frac{1}{\tilde{s}_r}}C_2.
\ea\ee
Note that  the right-hand side is independent of $m$ as well as $n$. 
Taking the limit as $i\to+\infty$, since $s_{i+1}\to +\infty$, by  the  convergence of the sums due to the d'Alambert criterion,  we conclude that
$$P_m\in L^\infty(Q_T) \ \mbox{for all } m\in \N\ \ \Rightarrow P \in L^\infty(Q_T).$$
From the relation $\pf=P+p_s$ and since $p_s\in L^\infty(Q_T),$ it finally follows that
$$\pf \in L^\infty(Q_T).$$

On the other hand, assuming \eqref{ass2} we first observe that if $p_s\in L^q(0,T; W^{1,q}(\O))$ with $q>10$ and $\bv \in L^{\frac{10}{3}}(Q_T)$, then $\bv\cdot \nabla p_s\in L^{\ell}(Q_T)$ with $\ell>\frac{5}{2}$. Consequently, \eqref{ass2} implies that $g+\bv\cdot \nabla p_s\in L^{\ell}(Q_T)$ with $\ell>\frac{5}{2}$. Then, we  conclude exactly as in the case given by \eqref{ass1} that 
$$\pf\in L^\infty(Q_T).$$
\chiu
The following lemma regards the attainment of the constitutive equations.
\begin{prop}[Convergence Lemma]\label{convergence-lemma}
Let $U\subset Q_T$ be an arbitrary measurable bounded set and let $\{\bZ^n\}_{n=1}^{+\infty}$, $\{\bD^n\}_{n=1}^{+\infty}$ and $\{\pf^n\}_{n=1}^{+\infty}$  
be such that 
\begin{align}\displ
&\bZ^n=\tau(\pf^n) \frac{\bD^n}{|\bD^n|+\frac{1}{n}} \mbox{ with } \tau(\pf^n)=q_*(p_s - \pf^n)^+,\label{C0}\\ \displ
&\sup_{n\in\N} \|\pf^n\|_\infty\leq C<\infty, \label{C1}\\\displ
&\bZ^n \rightharpoonup \bZ \ \textrm{weakly in } L^{2}(U)^{3\times 3},  \label{C2}\\\displ
&\bD^n \rightharpoonup \bD \ \textrm{weakly in } L^{2}(U)^{3\times 3}, \label{C3} \\\displ
&\pf^n  \to \pf \ \textrm{ strongly in }\  L^2(U) \mbox{ and a.e. in }U,\label{C4}\\\displ
&\limsup_{n\to \infty} \int_{U} \bZ^n : \bD^n \le \int_{U} \bZ : \bD,  \label{C5}
\end{align}
then, setting $\bS=\bZ+\bD$, 
\be\label{C6}
\bD =  \frac{\left(|\bS| - \tau(\pf)\right)^{+}}{|\bS|} \bS \ \ \mbox{a.e. in} \ U.
\ee
\end{prop}
\Pr We split the proof into three steps. Using the fact that \eqref{C6} is, by Proposition \ref{equivalent-formulation}, equivalent to \eqref{model2-2} (with $2\nu_*=q_*=1$), we first show that $|\bZ|\leq \tau(\pf). $ Then in order to verify that $\bZ:\bD\geq \tau(\pf)|\bD|$ in the third step, we show that $\bZ^n:\bD^n\rightharpoonup \bZ:\bD$ weakly in $L^1(U)$, which is the second part of the proof.
\\ \textit{Step 1.}  For all $n\in\N$,  by \eqref{C0}, $\bZ^n= \tau(\pf^n)\frac{ \mathbb{D}^n}{|\mathbb{D}^n|+\frac{1}{n}}$ and thus $|\bZ^n|\leq  \tau(\pf^n).$ For any subset $\omega\subset U$ it holds
\be\label{localrelation}
\int_\omega |\bZ^n|\leq \int_\omega \tau(\pf^n).
\ee
Since $\tau(\cdot)$ is Lipschitz, \eqref{C4} implies that
\be\label{CL1}
\tau(\pf^n)\to\tau(\pf) \ \mbox{ strongly in } L^2(U)
\ee
and
\be\label{CL1a}
\tau(\pf^n)\to\tau(\pf) \ \mbox{ a.e. in}\ U.
\ee
By virtue of \eqref{localrelation}, \eqref{CL1a} and the lower semicontinuity of $\int_\omega |\bZ^n|$ with respect to the weak convergence in $L^1(\omega)$ (which follows from \eqref{C2} since $U$ is bounded), we get
\be
\|\bZ\|_{L^1(\omega)}\leq \|\tau(\pf)\|_{L^1(\omega)} \ \ \mbox{for all  } \omega\subset U.
\ee
Lebesgue's Differentiation Theorem then implies
\be
|\bZ|\leq \tau(\pf) \  \mbox{a.e. in } U.
\ee
\textit{Step 2.} In order to establish that
\be \label{Pepa} 
\bZ^n:\bD^n\rightharpoonup \bZ:\bD \ \mbox{weakly in } L^1(U)
\ee
we set
\be\label{hat}
\widehat{\bZ}^n:= \tau(\pf^n) \frac{\bD}{|\bD|+\frac{1}{n}},
\ee
and
\be
\widehat{\bZ}:=\begin{cases} \tau(\pf) \frac{\bD}{|\bD|} &\mbox{if } \bD\neq\bO,\\\displ
\bO &\mbox{otherwise.} \end{cases}
\ee
\noindent
Thanks to \eqref{CL1a} we have that $\widehat{\bZ^n}\to \widehat{\bZ}$ almost everywhere in $Q_T$, and since $\widehat{\bZ}^n$ is  essentially bounded (because of \eqref{tauesb}), Lebesgue's Convergence Theorem yields
\begin{equation}\label{bZhL2}
\widehat{\bZ}^n\to \widehat{\bZ}\  \ \mbox{strongly in}\ L^2 (U).
\end{equation}
Employing   \eqref{C2} and \eqref{C5} and the convergences \eqref{bZhL2} and \eqref{C3}, we get
\begin{equation}\label{fin1}
 \limsup_{n\to \infty} \int_{U} (\bZ^n - \widehat{\bZ}^n):(\bD^n-\bD)\leq 0.
 \end{equation}
But since \eqref{hat}, the monotone property  \eqref{Prop-monotonicity} yields
$$(\bZ^n - \widehat{\bZ}^n):(\bD^n-\bD)\geq 0\ \mbox{a.e. in}\ U.$$
This together  with \eqref{fin1} implies that
$$(\bZ^n-\widehat{\bZ}^n):(\bD^n-\bD)\to 0 \ \mbox{strongly in}\ L^1(U),$$
and thus surely
\begin{equation}\label{fin2} (\bZ^n-\widehat{\bZ}^n):(\bD^n-\bD)\rightharpoonup 0\ \mbox{weakly in}\ L^1(U).\end{equation}
Since the strong convergence \eqref{bZhL2} and weak convergence \eqref{C3} imply that
\begin{equation}\label{finale}
\widehat{\bZ}^n:(\bD^n-\bD) \rightharpoonup 0\ \mbox{weakly in}\ L^1(U),
\end{equation}
so \eqref{fin2} yields
\begin{equation}\label{fin33}
\bZ^n:(\bD^n-\bD) \rightharpoonup 0\ \mbox{weakly in}\ L^1(U).
\end{equation}
Finally employing \eqref{C2} in \eqref{fin33} we conclude
\be
\bZ^n:\bD^n\rightharpoonup \bZ:\bD \ \mbox{weakly in } L^1(U).
\ee
\textit{Step 3.} It remains to show $\bZ:\bD\geq \tau(\pf)|\bD|$. First we note that
\be\label{difference}
|\tau(\pf^n)|\bD^n| - \bZ^n:\bD^n| = \tau(\pf^n) \frac{1}{n} \frac{ |\mathbb{D}^n|}{|\mathbb{D}^n|+\frac{1}{n}}.
\ee
Since $\tau$ is a Lipschitz function, \eqref{C1} gives
\be\label{tauesb} \|\tau(\pf^n)\|_\infty \leq C \ \mbox{uniformly in } \ n. \ee
Then the right hand side in \eqref{difference} is essentially bounded by $\frac{C}{n}$ and thus
\be\label{difference2}
|\tau(\pf^n)|\bD^n| - \bZ^n:\bD^n| \to 0 \ \mbox{in } L^\infty(U),
\ee
which implies that 
\be\label{mult-phi}
\lim_{n\to +\infty}\int_{U} \varphi (\tau(\pf^n)|\bD^n|-\bZ^n:\bD^n)=0 \ \ \mbox{for all } \varphi\in L^\infty(U).
\ee
Moreover, from  \eqref{C3} and \eqref{CL1} we get
\be
\varphi \tau(\pf^n) \bD^n\rightharpoonup \varphi\tau(\pf)\bD \ \mbox{in } L^1(U) \ \mbox{for all } \varphi\in L^\infty(U)
\ee
and the weak lower semicontinuity of the $L^1$-norm  implies that, for all $\varphi\in L^\infty(U)$ such that $\varphi\geq0$,
\be\label{liminf-phi}
\int_{U} \varphi\tau(\pf)|\bD| \leq \liminf_{n\to +\infty} \int_{U}\varphi \tau(\pf^n) |\bD^n|\,dx.  \ \  \ \ee
Using \eqref{Pepa} together with \eqref{mult-phi}, \eqref{liminf-phi} and \eqref{Pepa}, we obtain
\be
\int_{U} \varphi (\tau(\pf)|\bD|-\bZ:\bD)\leq \liminf_{n\to +\infty}\int_{U} \varphi (\tau(\pf^n)|\bD^n|-\bZ^n:\bD^n)=0
\ee
for any non-negative $\varphi\in L^\infty(U)$. Hence
$$
\bZ:\bD\geq \tau(\pf)|\bD|  \ \ \mbox{a.e. in } U, 
$$
which is \eqref{model2-2}$_2$.
\chiu

\section{Proof of the Main Theorem}\label{sectionProof}The proof is split in the following { five} steps.\\
\textit{Step 1. Approximations.} From Proposition \ref{prop} and Proposition \ref{prop52}, we get, for each $n\in \N$, the existence of $(\bv^n,\pf^n, \bZ^n, \bz^n)$ satisfying 
\be\label{533}\ba{l}\displ
\!\!\!\!\int_0^T\!\!\!\! \langle \partial_t \bv^n, \bw\rangle\,dt+\! \int_{Q_T} \!\!\!\! G_n\!\!\left(|\bv^n|^2\right) { \diver(\bv^n\!\otimes\!\bv^n \!) \cdot\! \bw)}\,dxdt+\!\!\! \int_{Q_T}\!\!\!\! (\bD \bv^n+\mathbb{Z}^{n}\!):\!\bD \bw \,dxdt 
\\\displ+\int_{\Sigma_T} \!\!\!\!(\bv_\tau^n+\bz^{n})\cdot  \bw_\tau\,d\sigma\,dt  -\int_0^T \!\!\!\langle\bb, \bw \rangle\,dt= 0
\mbox{ for all } \mathbf{w} \in L^2(0, T; W^{1,2}_{{\bf n},{\rm div}}),
\ea\ee
\be\label{534}\ba{l}\displ
\int_0^T\!\!\!\!\langle \partial_t \pf^{n}, z\rangle\,dt { -\!\!  \int_{Q_T} \!\!\!\! \pf^{n}\bv^{n}\cdot \nabla z}  + \nabla \pf^{n}\cdot \nabla z\,dx\,dt\! = \!\!\int_{Q_T}\!\!\! gz- { p_s \bv^{n}\!\! \cdot\!\! \nabla z}\,dxdt\\ \displ \hfill
\mbox{ for all } z\in L^2(0, T; W^{1,2}(\O)),
\ea\ee
and 
\begin{equation}\label{constitutive-eqs}
\bZ^n=(p_f^n-p_s)^+ \frac{\bD\bv^n}{|\bD\bv^n|+\frac{1}{n}} \textrm{ a.e. in  } Q_T \ \mbox{and } \bz^n= s_* \frac{\bv_{\tau}^n}{|\bv_{\tau}^n|+\frac{1}{n}} \mbox{ a.e. in } \Sigma_T.
\end{equation}
\\\textit{Step 2. Reconstruction of the pressure.} We set 
\begin{equation}\label{PR1}
p^n:= (-\Delta_{N})^{-1}\diver \mathbf{h}^n \ \mbox{ with } \int_\O p^n(t, \cdot)=0,
\end{equation}
where $-\Delta_{N}$ denotes the Laplace operator associated with the homogeneous Neumann boundary conditions and 
\begin{equation}\label{PR2}
\bh^n:= { - }\diver \left(  \bD \bv^n+ \bZ^n\right)+  \diver(\bv^n \otimes\bv^n) G_n(|\bv^n|^2) - \bb,
\end{equation}
associated with the boundary conditions $\bv\cdot\bn=0$ and $\bz^n=s_* \frac{\bv_{\tau}^n}{|\bv^n_\tau|+\frac{1}{n}}$ on $\Sigma_T$.
It means that $p^n$ solves (for a.a. $t\in [0,T]$)
\begin{equation}\label{PR3}
\begin{split}
(p^n, \Delta \varphi)& =\langle \bh^n, \nabla \varphi\rangle { + \int_{\partial \Omega} (\bv^n_{\tau} + \bz^n)\cdot (\nabla \varphi)_{\tau}} \\
&\qquad \mbox{ for all } \varphi\in W^{2,2}(\Omega) \ \mbox{ with } \nabla \varphi\cdot \bn=0 \mbox{ on } \partial \O,
\end{split}
\end{equation}
whereas 
\begin{equation}\label{PR4}
\bh^n\in L^2(0, T; \Wnm2).
\end{equation}
Consequently, 
\begin{equation}\label{PR5}
p^n\in L^2(0, T; L^2(\O)).
\end{equation}
Since any $\bw\in W^{1,2}_{\bn}$ satisfies 
$$
\bw = \tilde \bw + \nabla \varphi \, \textrm{ where }   \tilde{\bw}\in \Wnd2, \ \varphi\in W^{2,2}(\Omega), \nabla \varphi \cdot \bn = 0 \textrm{ on }\partial \Omega,
$$
we observe that due to \eqref{533} and \eqref{PR3} we get 
$$\ba{l}\displ (\bh^n, \bw)= (\bD\bv^n+\bZ^n, \bD\bw)+(\bv^n_\tau+\bz^n, \bw_\tau)_{\partial \Omega}+ (\diver(\bv^n\otimes \bv^n) G_n(|\bv^n|^2), \bw)- \langle\bb,  \bw \rangle\\
=- \langle \partial_t\bv^n, \tilde{\bw} \rangle+ (p^n,  \Delta \varphi)= - \langle \partial_t\bv^n, \tilde{\bw}+\nabla\varphi \rangle+ (p^n,  \diver(\tilde{\bw}+\nabla \varphi)).
\ea$$
Hence
\begin{equation}
\begin{split}
\langle\partial_t \bv^n, \bw\rangle &+(\bD \bv^n, \bD\bw)+ (\bZ^n, \bD\bw)  + (\bv^n , \bw_{\btau})_{\partial \Omega}+ (\bz^n, \bw_\tau)_{\partial \Omega}\\ \displ &+(\diver(\bv^n\otimes \bv^n) G_n(|\bv^n|^2), \bw)  
 = (p^n, \diver \bw) + \langle\bb,  \bw\rangle \, \textrm{ for all } \bw \in W^{1,2}_{\bn}.
 \end{split}\label{weak-n}
\end{equation}
\textit{Step 3. Uniform estimates with respect to $n$ and limit as $n\to+\infty$.}
Taking $\bv^n$ as test function in the  \eqref{533} and $\pf^n$ in the \eqref{534}, and proceeding similarly as in the derivation of \eqref{AEO} and \eqref{AE2} using also Korn's inequality, we obtain 
\begin{align}\label{unif-n}
&\sup_{n}\left(\|\bv^n\|_{L^\infty(L^2)} + \|\bD \bv^n\|_{L^2(Q_T)} + \|\nabla \bv^n\|_{L^2(Q_T)}+{ \|\bv^n_\tau\|}_{L^2(\Sigma_T)} \right) <+\infty,\\
\label{unif-pf}
&\sup_n \left(\|\pf^n\|_{L^\infty(0,T; L^2(\O))} + \|\nabla \pf^n \|_{L^2(Q_T)}\right) <+\infty,\\
&\sup_n\left( \|\bv^n\|_{\frac{10}{3}, Q_T}+\|\bv^n\|_{\frac{8}{3}, \Sigma_T}\right)<+\infty.\label{new0}
\end{align}
It follows directly from the proof of Proposition \ref{prop52} that
\begin{equation}\label{new00}
\sup_n\|\pf^n\|_{\infty, Q_T}<+\infty.
\end{equation}
It then follows from \eqref{constitutive-eqs} that 
\begin{equation}\label{newe1}
\sup_n\left(\|\bZ^n\|_{L^\infty(Q_T)}+\|\bz^n\|_{L^\infty(\Sigma_T)}\right)<+\infty.
\end{equation}
Since
$$ G_n(|\bv^n|^2)\diver (\bv^n\otimes \bv^n)=\sum_{h=1}^3 v^n_h\frac{\partial \bv^n}{\partial x_h}G_n(|\bv^n|^2),$$
and 
$\sup_n \|G_n(|\bv^n|^2)\|_{L^\infty(Q_T)}\leq 1$, it follows from \eqref{unif-n}, \eqref{new00} and H\"older's inequality that 
\begin{equation}\label{newe2}
\sup_n \|G_n(|\bv^n|^2)\diver (\bv^n\otimes \bv^n)\|_{L^{\frac{5}{4}}(Q_T)}<+\infty.
\end{equation}
For further analysis it is suitable to perform the following decomposition of the pressure $p^n$. Setting 
$$\bh^n_2:=G_n(|\bv^n|^2)\diver (\bv^n\otimes \bv^n)$$
and
$$p^n_2:= (-\Delta_{N})^{-1}\diver \mathbf{h}^n_2,$$
we conclude from \eqref{newe2} that 
\begin{equation}\label{newe3}
\sup_n \|\nabla p^n_2\|_{L^{\frac{5}{4}}(Q_T)}<+\infty.
\end{equation}
Furthermore, $\bh^n_1:=\bh^n-\bh^n_2$ fulfills $\sup_n\|\bh^n_1\|_{L^2(0,T; \Wnm2)}<+\infty,$ consequently $p^n_1:=p^n-p^n_2$ satisfies
\begin{equation}\label{newe4}
\sup_n \| p^n_1\|_{L^{2}(Q_T)}<+\infty,
\end{equation}
and it follows from \eqref{weak-n} that 
\be\label{omega}\ba{l}\displ\langle\partial_t \bv^n, \bw \rangle\!=\!\! \int_{\Omega}\!\! (-\bZ^n -\nabla \bv^n + p_1^n \bI)\!:\! \nabla \bw\,dx  
+ \langle \bb, \bw \rangle
 -\! \int_{\partial \Omega}\!\!\! (\bv^n_\tau+\bz^n)\! \cdot\! { \bw_{\tau}}\,d\sigma_x \\ \displ
 -\!\!\int_{\Omega}\!\! \!\left(G_n(|\bv^n|^2)\diver(\bv^n\! \otimes\! \bv^n) \!+\!\! \nabla p_2^n\right)\!\cdot\! \bw\,dx  \mbox{ for all } \bw\!\in\! L^2(0,T; \Wn2)\!\cap\! L^5(Q_T)^3\!.
\ea\ee
The above uniform estimates then imply that
\be \label{unif-vt}
\sup_n\|\partial_t v^n\|_{(L^2(0,T;W^{1,2}_{\bn}(\O))\cap L^{5}(Q_T)^3)^*}<+\infty,
\ee
and similarly
\be\label{pf-time}
\sup_n\|\partial_t \pf^n\|_{(L^2(0,T; W^{1,2}))^*}<+\infty.
\ee
Due to uniform estimates \eqref{unif-n}, \eqref{new00}, \eqref{newe1}, \eqref{newe2},  \eqref{newe3}, \eqref{newe4}, \eqref{unif-vt}, \eqref{pf-time}, the Aubin-Lions compactness lemma and the compact embedding of the Sobolev spaces into the space of traces, we get  the following convergences for subsequences that we do not relabel:
\begin{align}
&\bv^n \rightharpoonup \bv \textrm{ weakly in } L^2(0,T; W^{1,2}_{\bn}),\label{weak1}\\
&\pf^n \rightharpoonup \pf \textrm{ weakly in } L^2(0,T; W^{1,2}),\\
&\pf^n  \to \pf  \textrm{ strongly in } L^q(Q_T) \mbox{ for all } q\in\left[1, \frac{10}{3}\right),\\
&\partial_t\pf^n  \rightharpoonup \partial_t\pf \textrm{ weakly in } (L^2(0,T; W^{1,2}))^*,\\
&\pf^n \rightharpoonup^*  \pf \textrm{ weakly$^*$ in } L^\infty(Q_T),\\
&\bZ^n \rightharpoonup \bZ \textrm{ weakly$^*$ in } L^\infty(Q_T)^{3\times 3},\label{weak2}\\
&\bz^n \rightharpoonup \bz \textrm{ weakly$^*$ in } L^\infty(0,T; L^\infty(\partial \Omega)^3),\\
&\bv^n \to \bv \textrm{ a.e. in $Q_T$ and strongly in } L^q(Q_T)^3 \mbox{ for all } q\in\left[1, \frac{10}{3}\right),\label{strongL2}\\
&\bv^n_\tau \to \bv_\tau \textrm{ a.e. in $\Sigma_T$ and strongly in } L^q(0,T; L^q(\partial \Omega)^3)   \mbox{ for all } q\in\left[1, \frac{8}{3}\right),\label{vntau}\\
&p_1^n \rightharpoonup p_1 \textrm{ weakly in } L^2(Q_T),\\
&p_2^n \rightharpoonup p_2 \textrm{ weakly in } L^{\frac{5}{4}}(0,T; W^{1,\frac{5}{4}}(\Omega)),\label{unif-p2}\\
&\partial_t\bv^n \rightharpoonup \partial_t\bv \textrm{ weakly in }(L^2(0,T;W^{1,2}_{\bn}(\O))\cap L^{5}(Q_T)^3)^*,\\
&G_n(|\bv^n|^2)\diver (\bv^n\otimes \bv^n)\rightharpoonup \mathbf{g} \textrm{ weakly in } L^{\frac{5}{4}}(0,T; W^{1,\frac{5}{4}}(\Omega)).\label{div-unif}
\end{align}
It is not difficult to observe that due to the fact that
$$\|G_n(|\bv^n|^2)\|_{L^\infty(Q_T)}\leq 1 \mbox{ and } G_n(|\bv^n|^2)\to 1 \mbox{ strongly in } L^q(Q_T) \mbox{ for all } q\in [1,+\infty),$$
and due to \eqref{weak1} and \eqref{strongL2} we have
$$\mathbf{g}=\diver (\bv\otimes \bv).$$
Integrating then \eqref{weak-n} with respect the time between $0$ and $T$ and taking the limit as $n\to\infty$ we get
\be\label{a1}\ba{c}\displ
\!\!\!\!\int_0^T\!\!\!\! \!\langle \partial_t \bv, \bw\rangle\,dt+\!\! \int_{Q_T} \!\!\!\!\!\diver( \bv\otimes\bv)\!\! \cdot\! \bw\,dxdt+\!\!\! \int_{Q_T}\!\!\!\! \!\bD\bv\! :\!\bD \bw \,dxdt \\\displ
+\!\! \int_{Q_T} \!\!\!\!\mathbb{Z}\!:\! \bD\bw \, dxdt  -\!\!\int_0^T \!\!\!\!\langle\bb, \bw\rangle \,dt
+\! \int_{ \Sigma_T} \!\!\!\!(\bv_\tau+\bz)\!\cdot\!\bw_\tau\,d\sigma dt -\!\int_{Q_T}\!\! \!\!\!p_1\diver \bw\,dxdt
\\ \displ\hfill+\!\int_{Q_T}\!\! \!\!\!\nabla p_2\!\cdot\!\bw\,dxdt= 0 \mbox{ for all } \bw \in L^2(0,T;W^{1,2}_{\bn, \diver}(\O))\cap L^{5}(Q_T)^3.
\ea\ee
Integrating  \eqref{534}  with respect the time between $0$ and $T$ and taking the limit as $n\to\infty$ we get
\be\label{}\ba{c}\displ
\int_0^T\!\!\langle \partial_t \pf, z\rangle\,dt { -} \int_{Q_T} \!\!\!\!  \pf \bv\!\cdot\!\nabla z  +  \nabla \pf\!\cdot\! \nabla z\,dxdt = \int_{Q_T}\!\!\! gz\,dxdt
\\\displ
\hfill{ -}\!\! \int_{Q_T}\!\!\!\!  p_s\bv  \cdot \nabla z\,dxdt \mbox{   for all } z\in L^2(0,T; W^{1,2}(\O)).
\ea\ee
\textit{Step 4. Attainment of the constitutive equation on the boundary.} Since the structure of the constitutive equation on the boundary \eqref{constitutive-eqs}$_2$ is simpler than that used in Proposition \ref{convergence-lemma}, we can apply this assertion to this case as well. Indeed, we know that not only 
\begin{align*}
&\bz^n=s_*\frac{\bv^n_\tau}{|\bv^n_\tau|+\frac{1}{n}} \ \mbox{ a.e. on } \Sigma_T,\\
&\bz^n \rightharpoonup \bz  \mbox{ weakly in } L^q(\Sigma_T) \mbox{ for all } q\in[1, +\infty),\\
&\bv^n_\tau  \rightharpoonup \bv_\tau \textrm{ weakly in } L^{\frac{8}{3}}(\Sigma_T),
\end{align*}
but also
$$\bv^n_\tau \to \bv_\tau \mbox{ strongly in } L^q(\Sigma_T) \mbox{ for all } q\in\left[1, \frac{8}{3}\right).$$
Consequently,
$$ \lim_{n\to +\infty} \int_{\Sigma_T}\!\!\bz^n\!\cdot\!\bv^n_\tau d\sigma_xdt=\int_{\Sigma_T}\!\!\bz\!\cdot\!\bv_\tau d\sigma_xdt,$$
and by Proposition \ref{convergence-lemma} we get for $\bs=\bz+\bv_\tau$
$$\bv_\tau=\frac{(|\bs|-s_*)^+}{|\bs|}\bs \ \mbox{a.e. in } \Sigma_T.$$
\textit{Step 5. Attainment of the constitutive equation in the bulk.} We wish to use Proposition \ref{convergence-lemma}, and we can notice that all of its assumptions \eqref{C0}--\eqref{C4} are all fulfilled except \eqref{C5}. To prove it, we have to overcome the difficulty that $\bv$ is not an admissible test function in \eqref{a1}. This is why we employ the so-called $L^\infty$-truncation method applied to $\bv^n-\bv$.
Let $\{\lambda^n\}$, $A,B$  such that $0<A\leq \lambda^n\leq B<\infty$, where $A,B$ are independent of $n$ (but sufficiently large) and together with $\lambda^n$ will be specified later.
Consider the truncated velocity difference
\be
\bw^n :=T_{\lambda^n}(\bv^n - \bv) := (\bv^n - \bv) \min\left\{1, \frac{\lambda^n}{|\bv^n - \bv|} \right\}.
\ee
Since $\sup_n\|\bw^n\|_{\infty, Q_T}\leq B,$ and $\bv^n\to\bv$ a.e. in $Q_T$, Lebesgue's Theorem implies that
\be\label{strong-limit-wn}
\bw^n \to \bnul \ \mbox{strongly in } \ L^s(Q_T)^3 \ \mbox{for every } \ s\geq 1,
\ee
and similarly, as $\sup_n\|\bw^n\|_{\infty, Q_T}\leq B,$ and \eqref{vntau},
\be\label{wntau}\bw^n_\tau \to \bnul \ \mbox{ strongly in } \   L^2(0,T; L^2(\partial \Omega)^3).\ee
Since
\be\label{nablawn}
 \nabla\bw^n\!=\!\begin{cases}\!\nabla \bv^n-\nabla \bv \ \ \ \ \   \mbox{if} \ |\bv^n - \bv|<\lambda^n, \\ \displ
\!\!\frac{\lambda^n}{|\bv^n - \bv|}(\nabla \bv^n-\nabla \bv)-\!\lambda^n(\bv^n-\bv)\!\otimes\! \frac{(\nabla \bv^n-\nabla \bv)\! (\bv^n-\bv)}{|\bv^n-\bv|^3} &\mbox{otherwise},\end{cases}
\ee
we observe that 
\be\label{divwn}
|{\rm div} \ \bw^n| \leq \begin{cases}\displ 0 & \mbox{if} \ |\bv^n - \bv|<\lambda^n, \\ \displ
\frac{2 \lambda^n( |\nabla \bv^n| + |\nabla \bv|)}{|\bv^n - \bv|} & \mbox{otherwise}
\end{cases}
\ee
and 
\be\label{easy}
|\nabla\bw^n| \leq 2 |\nabla \bv^n-\nabla \bv|.
\ee
Then, due to \eqref{unif-n} and \eqref{easy}, $\nabla\bw^n$ is uniformly bounded in $L^2(Q_T)^{3\times 3}$ and, up to a subsequence, it converges weakly in $L^2(Q_T)^{3\times 3}$. But employing \eqref{strong-limit-wn} it follows that the weak limit has to be zero, i.e.
\be\label{weak-nablawn}
\nabla\bw^n \rightharpoonup \bO \ \mbox{ weakly in } L^2(Q_T)^{3\times 3} \mbox{ and } \bD\bw^n \rightharpoonup \bO \ \mbox{ weakly in } L^2(Q_T)^{3\times 3}.
\ee
Inserting $\bw^n$ in \eqref{omega}, we get
\be\label{limsup}\ba{l}\displ
\limsup_{n\to\infty}  \int_{Q_T} \!\!\!\!\mathbb{Z}^{n}\!:\!\! \nabla \bw^n - p^n_1 {\rm div}\ \bw^n+ \bD \bv^{n}\! :\!\bD \bw^n \,dxdt \\ \displ
=\limsup_{n\to\infty} \bigg[ \!-\!\!\int_0^T\!\!\!\! \langle \partial_t \bv^{n}, \bw^n\rangle\,dt- \int_{Q_T} \!\!\!\! G\left(|\bv^{n}|^2\right) {\rm div } \ (\bv^{n}\otimes\bv^{n})\!\cdot\!\bw^n\,dxdt  \\\displ
 +\int_{\Sigma_T} \!\!\!\!(\bv^n_\tau+\bz^{n})\!\cdot\!  \bw_\tau^n\,d\sigma dt  - \int_{Q_T}\!\!\! \nabla p^n_2\cdot \bw^n\,dxdt- \int_0^T \!\!\!\!\langle\bb, \bw^n\rangle dt\bigg].
\ea\ee
Now, by virtue of \eqref{strong-limit-wn}, \eqref{div-unif} and \eqref{unif-p2}, we observe that
\be\label{571}
\lim_{n\to\infty}  \int_{Q_T} \!\!\!\!\left(  G_n\left(|\bv^{n}|^2\right) {\rm div } (\bv^{n}\otimes\bv^{n}) \!+\! \nabla p^n_2 \right)\! \cdot\! \bw^n\,dxdt +\!\int_0^T\!\! \!\langle\bb, \bw^n\rangle dt= 0,
\ee
and by virtue of \eqref{unif-n} and \eqref{wntau}, it holds
\be\label{572}
\lim_{n\to\infty}  \int_{\Sigma_T}\!\!\! \! (\bv^n_\tau+\bz^{n})\!\cdot\!\bw_\tau^n d\sigma dt = 0.
\ee
Since $\bw^n\rightharpoonup \bnul$ weakly in $L^2(0,T; {W^{1,2}(\O)}^3)\cap {L^5(Q_T)}^3$ by \eqref{strong-limit-wn} then\\ $\displaystyle\lim_{n\to+\infty}  \int_0^T \langle\partial_t \bv, \bw^n\rangle\,dt=0$, thus
\be
\liminf_{n\to +\infty} \int_0^T \!\!\!\!  \langle\partial_t \bv^{n}, \bw^n\rangle\,dt= \liminf_{n\to +\infty} \int_0^T \!\!\! \langle  \partial_t (\bv^{n}-\bv),  \bw^n\rangle\,dt.
\ee
Moreover, 
\be\label{liminf}\ba{l}\displ
 \int_0^T \!\!\langle  \partial_t (\bv^{n}-\bv), \bw^n\rangle\,dt=  \int_0^T \!\! \partial_t\left(\! \frac{|\bv^{n}-\bv|^2}{2}\!\right)  \min\!\left\{1, \frac{\lambda^n}{|\bv^n - \bv|} \!\right\}\,dxdt\\\displ
\hfill =\int_{Q_T} \partial_t F^n(x,t)\,dx,
\ea\ee
where
$$ F^n(x,t):= \begin{cases}\displ
\frac{|\bv^{n}-\bv|^2}{2} & \mbox{if} \ |\bv^{n}-\bv|\leq \lambda^n, \\ \displ
\lambda^n |\bv^{n}-\bv|- \frac{{(\lambda^n)}^2}{2} & \mbox{if} \ |\bv^{n}-\bv|>\lambda^n.
\end{cases}$$
Thus from \eqref{liminf} and since $F^n(x,0)=0$ we get
\be
 \int_0^T \!\!\langle  \partial_t (\bv^{n}-\bv),\bw^n\rangle\,dt=\int_\Omega F^n(x,T)\,dx,
 \ee
taking the liminf we finally arrive at
\be\label{liminf2}
\liminf_{n\to +\infty} \int_0^T \!\!\!\langle \partial_t (\bv^{n}-\bv), \bw^n\rangle\,dt=\liminf_{n\to +\infty} \int_\Omega F^n(x,T)\,dx\geq0,
\ee
but this is equivalent to
\be\label{liminf3}
\limsup_{n\to\infty} \bigg[ \!-\!\!\int_0^T\!\!\!\! \langle \partial_t \bv^{n}, \bw^n\rangle\,dt\bigg]\leq 0.
\ee
Collecting  \eqref{571}, \eqref{572}, \eqref{liminf3}, it follows from \eqref{limsup} that
\be\label{limsup0}
\limsup_{n\to\infty}  \int_{Q_T} \!\!\!\!\mathbb{Z}^{n}\!:\!\bD \bw^n - \left(p^n_1 {\rm div}\ \bw^n \right)+ \bD \bv^{n}\! :\!\bD \bw^n\,dxdt\leq 0.
\ee
Since $(\bD\bv^n-\bD\bv)\!:\!\bD\bw^n\geq 0$ and $\lim_{n\to+\infty} \int_{Q_T}\!\!\!\!\bD\bv\!:\!\bD\bw^ndxdt=0$, \eqref{divwn} and \eqref{limsup0} imply that 
\be\label{limsup1}\begin{split}
\limsup_{n\to\infty}  \int_{Q_T} \!\!\!\!\mathbb{Z}^{n}\!:\!\bD \bw^n &+ \bD \bv^{n}\! :\!\bD \bw^n\,dxdt \\
&\leq \limsup_{n\to\infty}  \int_{Q_T} \!\!  |p^n_1| |{\rm div}\ \bw^n| \,dxdt \\
&\leq  \limsup_{n\to\infty}\int_{\{|\bv^n-\bv|\geq \lambda^n\}}\frac{\lambda^n}{|\bv^n-\bv|}|p_1^n|(|\nabla \bv^n|+|\nabla \bv|).
\end{split}
\ee
Let  $\overline{\bZ}\in {L^{\frac{10}{3}}(Q_T)}^{3\times 3}$  be such that
\be\label{Zbar}\overline{\bZ}=
\begin{cases}\vspace{6pt}\displaystyle
\bO \ &\mbox{if} \  \bD\bv=\bO ,\\\vspace{6pt}\displaystyle
\tau(\pf)\frac{ \bD}{|\bD|}  &\mbox{if} \ \bD\bv\neq\bO. \end{cases}
\ee
Since $\lim_{n\to+\infty} \int_{Q_T}\overline{\bZ}\!:\!\bD\bw^ndxdt=0$ thanks to \eqref{weak-nablawn}, we arrive at 
\begin{equation}\label{new}
\limsup_{n\to\infty}  \int_{Q_T} \!\!\!\!(\mathbb{Z}^{n}-\overline{\bZ})\!:\!\bD \bw^n \,dxdt\leq \limsup_{n\to\infty}  \int_{\{|\bv^n-\bv|\geq \lambda^n\}}\frac{\lambda^n}{|\bv^n-\bv|}|p_1^n|(|\nabla \bv^n|+|\nabla \bv|).
\end{equation}
Splitting the integral on the left-hand side of \eqref{new}  into  two parts, one integrated over $\{|\bv^n-\bv|\leq \lambda^n\}$ the other over $\{|\bv^n-\bv|\geq \lambda^n\}$, using \eqref{nablawn}, and moving the latter to the right-hand side and estimating it by \eqref{easy}, we get
\be\label{584}
\limsup_{n\to \infty}\int_{\{|\bv^n-\bv|\leq \lambda^n\}}\!\!\! (\bZ^n-\overline{\bZ}) : \bD(\bv^n - \bv)
\leq C \limsup_{n\to\infty}\int_{\{|\bv^n-\bv|\geq \lambda^n\}}\frac{\lambda^n}{|\bv^n-\bv|}I^n,\ee
where
$$
  I^n:= (|p_1^n|^2 + |\nabla \bv^n|^2 +|\nabla \bv|^2 + |\bZ^n|^2+|\overline{\bZ}|^2 ) \mbox{ and } \sup_n\int_{Q_T}\!\!\!\!\! I^n \,dxdt<+\infty.
$$
Let $N\in \mathbb{N}$ be arbitrary. We fix $A=N$ and $B=N^{N+1}$ and define
\begin{equation*}
Q_i^n:= \{(t,x)\in Q_T; N^i \le |\bv^n - \bv| \le N^{i+1}\}, \qquad i=1,\dots, N.
\end{equation*}
Since
\begin{equation}
\sum_{i=1}^N \int_{Q_i^n} I^n \le C_*,
\end{equation}
there is, for each $n\in \mathbb{N}$, an index $i_n\in\{1, \dots, N\}$ such that
\begin{equation}\label{586}
\int_{Q_{i_n}^n} I^n \le \frac{C_*}{N}.
\end{equation}
Setting $\lambda^n=N^{i_n}$, the right-hand side of \eqref{584} can be estimated as follows using  \eqref{586} and the fact that $I^n$ is uniformly bounded in $L^1(Q_T)$
\begin{equation}\label{587}
\begin{split}
\int_{\{|\bv^n-\bv|\ge N^{i_n}\}} & \frac{N^{i_n}}{|\bv^n-\bv|} I^n  \\ = &\int_{\{N^{i_n}\le |\bv^n-\bv|\le N^{i_n+1}\}} \frac{N^{i_n}}{|\bv^n-\bv|} I^n 
 + \int_{\{|\bv^n-\bv|\ge N^{i_n+1}\}} \frac{N^{i_n}}{|\bv^n-\bv|} I^n  \\
\le &\int_{Q_{i_n}^n} I^n + \frac{1}{N} \int_{\{|\bv^n-\bv|\ge N^{i_n+1}\}} I^n \le \frac{C_*}{N}.
\end{split}
\end{equation}
Let 
$$
W^n:=  \left(\bZ^n - \overline{\bZ} \right)\!:\! \left( \bD\bv^n-\bD\bv\right).
$$
Then   \eqref{584} and \eqref{587} imply that 
\begin{equation}\label{588}\ba{c}\displ
\limsup_{n\to \infty}\! \int_{|\bv^n-\bv|\le \lambda^n}\!\!\!\! \bW^n \le \frac{C_*}{N}\\\displ
\iff \limsup_{n\to \infty}\! \int_{|\bv^n-\bv|\le \lambda^n}\!\!\!\! |\bW^n| \le \frac{C_*}{N} +  2\limsup_{n\to \infty}\! \int_{|\bv^n-\bv|\le \lambda^n}\!\!\!\! {(\bW^n)}^-
\ea\end{equation}
Now we show that 
\be\label{589}
{(\bW^n)}^-\to 0 \ \mbox{strongly in } \ L^1(Q_T).
\ee
Recalling that $\bZ^n=\mathcal{Z}^n( \pf^n, \bD\bv^n)$ and incorporating  \eqref{Prop-monotonicity}, we get 
\be
\begin{split}
\bW^n= &(\bZ^n-\mathcal{Z}^n( \pf^n, \bD\bv))\! :\! (\bD\bv^n-\bD\bv) + (\mathcal{Z}^n( \pf^n, \bD\bv)-\overline{\bZ})\!:\!(\bD\bv^n-\bD\bv) \\&
\geq (\mathcal{Z}^n(\bD\bv, \tau({\pf^n}))-\overline{\bZ})\! :\! (\bD\bv^n-\bD\bv).
\end{split}
\ee
Splitting $Q_T=\{ |\bD\bv|=0\}\cup\{ |\bD\bv|>0\}$, thanks to the definitions of $\bZ^n$ and $\overline{\bZ}$ and since $\pf^n$ converges pointwise, we get
$$\mathcal{Z}^n(\bD\bv, \tau({\pf^n}))\to \overline{\bZ} \ \mbox{a.e. in } 	\ Q_T. $$
Also, independently of $n$,
$$ |\mathcal{Z}^n(\bD\bv, \tau({\pf^n})) - \overline{\bZ}|\leq C|\bD\bv|.$$
By the Dominated Convergence Theorem and since $(\bD\bv^n-\bD\bv)$ is bounded in $L^2(Q_T)$ yield \eqref{589}. \\
Combining  \eqref{588}, \eqref{589} and recalling that $A=N\leq \lambda^n$,
\begin{equation}\label{591}
\limsup_{n\to \infty} \int_{|\bv^n-\bv|\le N}  |\bW^n| \le \frac{C_*}{N}\,.
\end{equation}
With the help of the H\"older and Chebyshev inequalities, we observe that 
\be\ba{l}\displ
\int_{Q_T} \sqrt{|\bW^n|} \leq  \int_{|\bv^n-\bv|\leq N} \sqrt{|\bW^n|} + \int_{|\bv^n-\bv|> N}\sqrt{|\bW^n|} \\\displ
\leq |Q_T|^{\frac{1}{2}} \sqrt{ {\int_{|\bv^n-\bv|\le N} |\bW^n|} }+ \|\bW^n\|_{L^2(Q_T)}^{\frac{1}{2}} \sqrt{ |\{ |\bv^n-\bv|> N\}| }\leq \frac{C}{\sqrt{N}}
\ea
\ee
which implies that for a suitable subsequence,
\be
\bW^n \to 0 \quad \textrm{ a.e. in } Q_T\,.
\ee
Applying  Egorov Theorem, one concludes that
\begin{equation*}
W^n \to 0 \quad \textrm{ strongly in } L^1(Q_T\setminus E_j)\,,
\end{equation*}
where $E_j\subset Q_T$ are such that $\lim_{j\to \infty} |E_j| = 0$.
It follows from the definition of $W^n$ and the weak convergences \eqref{weak1}, \eqref{weak2} that
\begin{equation*}
\limsup_{n\to \infty} \int_{Q_T\setminus E_j} \!\!\!\bZ^n \!:\! \bD\bv^ndxdt
= \limsup_{n\to \infty} \int_{Q_T\setminus E_j}\!\!\!\!\!\! \overline{\bZ}\!:\! (\bD\bv^n-\bD\bv) + \bZ^n\!:\!\bD\bv dxdt =  \int_{Q_T\setminus E_j}\!\!\!\!\!\!\!\!\! \bZ \!:\! \bD\bv \,dxdt.
\end{equation*}
Thus, the assumptions \eqref{C0}-\eqref{C5} of Proposition \ref{convergence-lemma} are verified with $U=Q_T\setminus E_j$, for all $j\in \N$. Due to the properties of $E_j$, we finally conclude, using \eqref{C6}, that
$$\bD\bv=\frac{(|\bS|-\tau(\pf))^+}{|\bS|}\bS.$$
The proof of Theorem \ref{main-thm} is complete.
\section*{Conclusion}
This study has been inspired by recent research concerning implicitly constituted materials on one hand and by a recent interesting paper by Chupin and Math\'e \cite{CM} on the other hand. This study extends the results presented in \cite{CM} in several directions. First, we have studied slightly different system of PDEs, namely the one we were able to derive from the basic governing equations of the theory of mixtures, under  the cascade of several justified simplifications. Second, the activated system contains in comparison to \cite{CM}, a non-trivial right-hand side in the equation for the fluid pressure $\pf$. Consequently, we had  to use a different approach to get $L^\infty$-estimates for $\pf$. Third, inspired by \cite{CM} we provide characterization of the constitutive equation in Proposition \ref{equivalent-formulation}. Using one of these equivalent descriptions, one can correct the proof in \cite{CM} and get a useful tool exploited in the proof of Proposition \ref{convergence-lemma} here. Fourth, we considered stick-slip boundary conditions that are not only physically relevant but, on contrary to no-slip boundary condition, guarantees the integrability of the pressure up to the boundary. Finally, we use $L^\infty$-truncation method to analyze three-dimensional flows (while the result in \cite{CM} concerns planar flows). 
We wish to remark that it is possible to use, instead of $L^\infty-truncation$, a solenoidal Lipschitz truncation (introduced in \cite{BreitDienSc}) and consider the formulation free of the pressure. 
From the application point of view, the system of equations analyzed here { has some relevance to the problem of static liquefaction and enhanced oil recovery}. Of course, as one may conclude from Section 2, this topic provides several questions for further research.
\section*{Acknowledgements}
The authors thank Miroslav Bul{\'\i}{\v{c}}ek for several valuable comments and discussions.
A.~Abbatiello is partially supported by the Italian National Group of Mathematical Physics (GNFM-INdAM)  via GNFM Progetto Giovani 2017 and, is also grateful to Charles University for the hospitality during her stay when the work was performed.
T.~Los, J.~M\'alek,  and O.~Sou\v{c}ek acknowledge support of the project 18-12719S financed by the Czech Science Foundation. T.~Los is also thankful to the institutional support  through the project GAUK 550218.
\nocite{*}
\providecommand{\bysame}{\leavevmode\hbox to3em{\hrulefill}\thinspace}
\providecommand{\MR}{\relax\ifhmode\unskip\space\fi MR }
\providecommand{\MRhref}[2]{%
  \href{http://www.ams.org/mathscinet-getitem?mr=#1}{#2}
}
\providecommand{\href}[2]{#2}

\end{document}